\documentclass[12pt]{amsart}
\usepackage{amssymb}

\usepackage{amscd}
\newtheorem{theorem}{Theorem}[section]

\newtheorem{lemma}[theorem]{Lemma}

\newtheorem{remark}{Remark}[section]

\allowdisplaybreaks
\theoremstyle{definition}
\theoremstyle{remark}
\numberwithin{equation}{section}

\let\que=q
\begin{document}
\title[{\small $\que$}-Fourier Series]{Basic Analog of Fourier Series\\
on a {\large $\que$}-quadratic grid}
\author{Joaquin Bustoz}
\address{\hskip-\parindent
Joaquin Bustoz\\
Department of Mathematics\\
Arizona State University\\
Tempe, Arizona 85287-1804, U.S.A.}
\email{bustoz@math.la.asu.edu}
\author{Sergei K. Suslov}
\address{\hskip-\parindent 
Sergei K. Suslov\\
Kurchatov Institute\\
Moscow, 123182, Russia}
\curraddr{Department of Mathematics\\
Arizona State University \\
Tempe, Arizona 85287-1804\\
U.S.A.}
\email{suslov@math.la.asu.edu}
\date{May 21, 1997}
\subjclass{Primary 33B10, 33D15; Secondary 42C10}
\keywords{Trigonometric functions, basic trigonometric functions, orthogonality
relations, Fourier series, $q$-Fourier series}

\begin{abstract}
We prove orthogonality relations for some analogs of trigonometric functions
on a $q$-quadratic grid and introduce the corresponding $q$-Fourier series.
We also discuss several other properties of this basic trigonometric system
and the $q$-Fourier series.
\end{abstract}

\maketitle
\section{Introduction}

A periodic function with period $2l$,

\begin{equation}
f(x+2l)=f(x),
\end{equation}
can be represented as the Fourier series,

\begin{equation}
f(x)=a_{0}+\sum_{n=1}^{\infty }\left( a_{n}\cos \frac{\pi n}{l}x+b_{n}\sin 
\frac{\pi n}{l}x\right) ,
\end{equation}
where 
\begin{eqnarray}
&&a_{0}=\frac{1}{2l}\ \int_{-l}^{l}f(x)\ dx, \\
&&a_{n}=\frac{1}{l}\ \int_{-l}^{l}f(x)\cos \frac{\pi n}{l}x\ dx, \\
&&b_{n}=\frac{1}{l}\ \int_{-l}^{l}f(x)\sin \frac{\pi n}{l}x\ dx.
\end{eqnarray}
For convergence conditions of (1.2) see, for example, \cite{AKh}, \cite
{Wh:Wa}, and \cite{Zyg}. The formulas (1.3)--(1.5) for the coefficients of
the Fourier series are consequences of the orthogonality relations for
trigonometric functions 
\begin{eqnarray}
&&\int_{-l}^{l}\cos \frac{n\pi x}{l}\ \cos \frac{m\pi x}{l}\ dx=0,\qquad \
m\neq n, \\
&&\int_{-l}^{l}\sin \frac{n\pi x}{l}\ \sin \frac{m\pi x}{l}\ dx=0,\qquad \
m\neq n, \\
&&\int_{-l}^{l}\cos \frac{n\pi x}{l}\ \sin \frac{m\pi x}{l}\ dx=0,\qquad \
m\neq n.
\end{eqnarray}
In the present paper we discuss a $q$-version of the Fourier series (1.1)
with the aid of basic or $q$-analogs of trigonometric functions introduced
recently in \cite{Is:Zh} (see also \cite{At:Su} and \cite{Su}). Our first
main objective will be to establish analogs of the orthogonality relations
(1.5)--(1.7) for $q$-trigonometric functions on a $q$-quadratic grid.

There are several ways to prove the orthogonality relations (1.6)--(1.8) for
trigonometric functions. The method based on the second order differential
equation, 
\begin{equation}
u^{\prime \prime }+\omega ^{2}u=0,
\end{equation}
can be extended to the case of basic trigonometric functions. Consider, for
example, two functions $\cos \omega x$ and $\cos \omega ^{\prime }x$, which
satisfy (1.9) with different eigenvalues $\omega $ and $\omega ^{\prime }$.
Then, 
\begin{eqnarray}
&\left( \omega ^{2}-\omega ^{\prime 2}\right) &\!\int_{-l}^{l}\cos \omega x\
\cos \omega ^{\prime }x\ dx \\
&&\ \qquad \qquad =W\left( \cos \omega x,\ \cos \omega ^{\prime }x\right)
|_{-l}^{l}  \notag \\
&&\qquad \qquad \quad =\left. \left| 
\begin{array}{ll}
\ \ \cos \omega x & \ \ \cos \omega ^{\prime }x \\ 
-\omega \ \sin \omega x & -\omega ^{\prime }\ \sin \omega ^{\prime }x
\end{array}
\right| \ \right| _{-l}^{l}.  \notag
\end{eqnarray}
The right side of (1.10) vanishes when 
\begin{equation}
\sin \omega l=\sin \omega ^{\prime }l=0,
\end{equation}
which gives 
\begin{equation}
\omega =\frac{\pi }{l}\ n,\qquad \omega ^{\prime }=\frac{\pi }{l}\ m,
\end{equation}
where $n,m=0,\pm 1,\pm 2,\pm 3,...\ $. In the same manner, one can prove
(1.7). The last equation (1.8) is valid by symmetry. We shall extend this
consideration to the case of the basic trigonometric functions in the
present paper.

This paper is organized as follows. In Section 2 we introduce the $q$%
-trigonometric functions. In the next section we derive a continuous
orthogonality property of these functions, and then, in Section 4, we
formally discuss the limit $q\rightarrow 1$ of these new orthogonality
relations. Section 5 is devoted to the investigation of some properties of
zeros of the basic trigonometric functions and in Section 6 we evaluate the
normalization constants in the orthogonality relations for these functions.
In Section 7 we state the orthogonality relation for the corresponding $q$%
-exponential functions. Finally, we introduce basic analogs of Fourier
series in Section 8, and in Sections 9--11 we give a proof of the
completeness of the $q$-trigonometric system and establish some elementary
facts about convergence of our $q$-Fourier series. Examples of these series
are considered in Sections 12 and 14; we prove some basic trigonometric
identities we heed in Section 13. Some miscellaneous results concerning $q$%
-trigonometric functions are discussed in Section 15.

\section{Analogs of Trigonometric Functions on a $q$-Quadratic Grid}

The following functions $C(x)$ and $S(x)$ given by 
\begin{eqnarray}
C(x) &=&C_{q}(x;\omega ) \\
&=&\frac{\left( -\omega ^{2};q^{2}\right) _{\infty }}{\left( -q\omega
^{2};q^{2}\right) _{\infty }}\   \notag \\
&&\times \ _{2}\varphi _{1}\left( 
\begin{array}{c}
\begin{array}{ll}
-qe^{2i\theta }, & -qe^{-2i\theta }
\end{array}
\\ 
q
\end{array}
;\ q^{2},\ -\omega ^{2}\right)  \notag
\end{eqnarray}
and 
\begin{eqnarray}
S(x) &=&S_{q}(x;\omega ) \\
&=&\frac{\left( -\omega ^{2};q^{2}\right) _{\infty }}{\left( -q\omega
^{2};q^{2}\right) _{\infty }}\ \frac{2q^{1/4}\omega }{1-q}\ \cos \theta 
\notag \\
&&\times \ \ _{2}\varphi _{1}\left( 
\begin{array}{c}
\begin{array}{ll}
-q^{2}e^{2i\theta }, & -q^{2}e^{-2i\theta }
\end{array}
\\ 
q^{3}
\end{array}
;\ q^{2},\ -\omega ^{2}\right) ,  \notag
\end{eqnarray}
were discussed recently \cite{At:Su}, \cite{Is:Zh} and \cite{Su} as $q$%
-analogs of $\cos \omega x$ and $\sin \omega x$ on a $q$-quadratic lattice $%
x=\cos \theta $.

These functions are special cases $y=0$ of more general basic trigonometric
functions 
\begin{eqnarray}
&&C(x,y)=C_{q}(x,y;\omega )  \notag \\
&&\ \!\!=\frac{\left( -\omega ^{2};q^{2}\right) _{\infty }}{\left( -q\omega
^{2};q^{2}\right) _{\infty }}\   \notag \\
&&\times \ _{4}\varphi _{3}\left( 
\begin{array}{c}
\!-q^{1/2}e^{i\theta +i\varphi },\!-q^{1/2}e^{i\theta -i\varphi
},\!-q^{1/2}e^{i\varphi -i\theta },\!-q^{1/2}e^{-i\theta -i\varphi } \\ 
-q,\quad q^{1/2},\quad -q^{1/2}\quad
\end{array}
;q,\!-\omega ^{2}\right)  \notag \\
&&\qquad
\end{eqnarray}
and 
\begin{eqnarray}
&&S(x,y)=S_{q}(x,y;\omega )  \notag \\
&&\ =\frac{\left( -\omega ^{2};q^{2}\right) _{\infty }}{\left( -q\omega
^{2};q^{2}\right) _{\infty }}\ \frac{2q^{1/4}\omega }{1-q}\ \left( \cos
\theta +\cos \varphi \right)  \notag \\
&&\ \times \ _{4}\varphi _{3}\left( 
\begin{array}{c}
-qe^{i\theta +i\varphi },-qe^{i\theta -i\varphi },-qe^{i\varphi -i\theta
},-q^{-i\theta -i\varphi } \\ 
-q,\quad q^{3/2},\quad -q^{3/2}
\end{array}
;\ q,\ -\omega ^{2}\right)  \notag \\
&&\qquad
\end{eqnarray}
which are $q$-analogs of $\cos \omega (x+y)$ and $\sin \omega (x+y)$,
respectively (see \cite{Su}). Here $x=\cos \theta $ and $y=\cos \varphi $.
Usually we shall drop $q$ from the symbols $C_{q}(x;\omega )$, $%
S_{q}(x;\omega )$, $C_{q}(x,y;\omega )$, and $S_{q}(x,y;\omega )$ because
the same base is used throughout the paper.

The symbols $_{2}\varphi _{1}$ and $_{4}\varphi _{3}$ in (2.1)--(2.4) are,
of course, special cases of basic hypergeometric functions, 
\begin{eqnarray}
&&\ _{r}\varphi _{s}\left( 
\begin{array}{c}
a_{1},\ a_{2},...\ ,\ a_{r} \\ 
b_{1},\ b_{2},...\ ,\ b_{s}
\end{array}
;\ q,\ t\right) \\
&&\quad :=\sum_{n=0}^{\infty }\frac{\left( a_{1},\ a_{2},...\ ,\ a_{r};\
q\right) _{n}}{\left( q,\ b_{1},\ b_{2},...\ ,\ b_{s};\ q\right) _{n}}\
\left( \left( -1\right) ^{n}q^{n\left( n-1\right) /2}\right) ^{1+s-r}\ t^{n}.
\notag
\end{eqnarray}
The standard notations for the $q$-shifted factorial are 
\begin{equation}
\left( a;q\right) _{0}:=1,\qquad \left( a;q\right)
_{n}:=\prod_{k=0}^{n-1}\left( 1-aq^{k}\right) ,
\end{equation}
\begin{equation}
\left( a_{1},\ a_{2},...,\ a_{m};q\right) _{n}:=\prod_{l=1}^{m}\left(
a_{l};q\right) _{n},
\end{equation}
where $n=1,2,...,$ or $\infty $, when $|q|<1$. See \cite{Ga:Ra} for an
excellent account of the theory of basic hypergeometric functions.

Functions (2.1)--(2.4) are defined here for $|\omega |<1$ only. For an
analytic continuation of these functions in a larger domain see \cite
{Is:Ma:Su1}, \cite{Is:Zh}, and \cite{Su}. For example, 
\begin{eqnarray}
C(x) &=&\frac{\left( q\omega ^{2}e^{2i\theta },\ q\omega ^{2}e^{-2i\theta
};\ q^{2}\right) _{\infty }}{\left( q,\ -q\omega ^{2};\ q^{2}\right)
_{\infty }} \\
&&\times \ _{2}\varphi _{2}\left( 
\begin{array}{cc}
-\omega ^{2}, & -q\omega ^{2} \\ 
q\omega ^{2}e^{2i\theta }, & q\omega ^{2}e^{-2i\theta }
\end{array}
;\ q^{2},\ q\right)  \notag
\end{eqnarray}
and 
\begin{eqnarray}
S(x) &=&\frac{\left( q^{2}\omega ^{2}e^{2i\theta },\ q^{2}\omega
^{2}e^{-2i\theta };\ q^{2}\right) _{\infty }}{\left( q^{3},\ -q\omega ^{2};\
q^{2}\right) _{\infty }}\ \frac{2q^{1/4}\omega }{1-q}\ \cos \theta \\
&&\times \ _{2}\varphi _{2}\left( 
\begin{array}{cc}
-\omega ^{2}, & -q\omega ^{2} \\ 
q^{2}\omega ^{2}e^{2i\theta }, & q^{2}\omega ^{2}e^{-2i\theta }
\end{array}
;\ q^{2},\ q^{3}\right) .  \notag
\end{eqnarray}
One can see from (2.8) and (2.9) that the basic trigonometric functions
(2.1) and (2.2) are entire functions in $z$ when $e^{i\theta }=q^{z}$.
Analytic continuation of $q$-trigonometric functions (2.3) and (2.4) can be
obtained on the basis of the ``addition'' theorems, 
\begin{eqnarray*}
C\left( x,y\right) &=&C\left( x\right) C\left( y\right) -S\left( x\right)
S\left( y\right) , \\
S\left( x,y\right) &=&S\left( x\right) C\left( y\right) +C\left( x\right)
S\left( y\right) ,
\end{eqnarray*}
found in \cite{Su}.

The basic trigonometric functions (2.1)--(2.4) are solutions of a difference
analog of equation (1.9) on a $q$-quadratic lattice, 
\begin{equation}
\sigma \ \frac{\Delta }{\nabla x_{1}(z)}\left( \frac{\nabla u(z)}{\nabla x(z)%
}\right) +\lambda \ u(z)=0,
\end{equation}
where $x(z)=\frac{1}{2}\left( q^{z}+q^{-z}\right) $, $q^{z}=e^{i\theta }$, $%
x_{1}(z)=x(z+1/2)$, $\lambda /\sigma =4q^{1/2}\omega ^{2}/\left( 1-q\right)
^{2}$, and $\Delta f(z)=\nabla f(z+1)=f(z+1)-f(z)$. See \cite{At:Su}, \cite
{Is:Zh}, \cite{Ni:Su:Uv}, \cite{Su0}, and \cite{Su} for more details.
Equation (2.10) can also be rewritten in a more symmetric form, 
\begin{equation}
\ \frac{\delta }{\delta x(z)}\left( \frac{\delta u(z)}{\delta x(z)}\right) +%
\frac{4q^{1/2}}{\left( 1-q\right) ^{2}}\ \omega ^{2}\ u(z)=0,
\end{equation}
where $\delta f(z)=f(z+1/2)-f(z-1/2)$.

The $q$-trigonometric functions (2.1)--(2.4) satisfy the
difference-differentiation formulas 
\begin{equation}
\frac{\delta }{\delta x}\ C(x,y)=-\frac{2q^{1/4}}{1-q}\ \omega \ S(x,y)
\end{equation}
and 
\begin{equation}
\frac{\delta }{\delta x}\ S(x,y)=\frac{2q^{1/4}}{1-q}\ \omega \ C(x,y).
\end{equation}
See \cite{Is:Zh} and \cite{Su}. Applying the operator $\delta /\delta x$ to
the both sides of (2.12) or (2.13) we obtain equation (2.11) again.

Equation (2.10) is a very special case of a general difference equation of
hypergeometric type on nonuniform lattices (cf. \cite{At:Su}, \cite{Ni:Su:Uv}%
, \cite{Su0}, and \cite{Su}). The Askey--Wilson polynomials and their
special and limiting cases \cite{As:Wi}, \cite{Ko:Sw}, and \cite{Ni:Su:Uv}
are well-known as the simplest and the most important orthogonal solutions
of this difference equation of hypergeometric type. Recently, Ismail,
Masson, and Suslov \cite{Is:Ma:Su1}, \cite{Is:Ma:Su2}, \cite{Su1}, \cite{Su2}
have found another type of orthogonal solutions of this difference equation.
In the present paper we shall discuss this new orthogonality property at the
level of basic trigonometric functions.

\section{Continuous Orthogonality Property for $q$-Trigonometric Functions}

Our main objective in this paper is to find the orthogonality relations for $%
q$-trigonometric functions (2.1)--(2.2) similar to the orthogonality
relations (1.5)--(1.7). Consider difference equations for the functions $%
u(z)=C_{q}\left( x(z);\omega \right) $ and $v(z)=C_{q}\left( x(z);\omega
^{\prime }\right) $ in self-adjoint form, 
\begin{equation}
\frac{\Delta }{\nabla x_{1}(z)}\left( \sigma \rho (z)\ \frac{\nabla u(z)}{%
\nabla x(z)}\right) +\lambda \ \rho (z)u(z)=0
\end{equation}
and 
\begin{equation}
\frac{\Delta }{\nabla x_{1}(z)}\left( \sigma \rho (z)\ \frac{\nabla v(z)}{%
\nabla x(z)}\right) +\lambda ^{\prime }\ \rho (z)v(z)=0,
\end{equation}
where the function $\rho (z)$ satisfies the ``Pearson equation''\cite
{Ni:Su:Uv}, \cite{Su0}, 
\begin{equation}
\frac{\rho (z+1)}{\rho (z)}=\frac{\sigma (-z)}{\sigma (z+1)}=1=q^{-4z-2}\
\left( q^{2}\right) ^{2z+1},
\end{equation}
and 
\begin{equation}
\lambda =\frac{4q^{1/2}\sigma }{\left( 1-q\right) ^{2}}\ \omega ^{2},\qquad
\lambda ^{\prime }=\frac{4q^{1/2}\sigma }{\left( 1-q\right) ^{2}}\ \omega
^{\prime 2}.
\end{equation}

One can easily check that 
\begin{equation}
\frac{\rho _{0}(z+1)}{\rho _{0}(z)}=q^{-4z-2}\quad \text{for\quad }\rho
_{0}(z)=\frac{\left( q^{2z},\ q^{-2z};\ q\right) _{\infty }}{q^{z}-q^{-z}}
\end{equation}
and 
\begin{equation}
\frac{\rho _{\alpha }(z+1)}{\rho _{\alpha }(z)}=q^{4z+2}\quad \text{for\quad 
}\rho _{\alpha }(z)=\left( q^{2\alpha +2z},\ q^{2\alpha -2z},\ q^{2-2\alpha
+2z},\ q^{2-2\alpha -2z};\ q^{2}\right) _{\infty }^{-1}
\end{equation}
(cf. \cite{Is:Ma:Su1}, \cite{Su1}, and \cite{Su2}). Therefore, we can choose
the following solution of (3.3), 
\begin{equation}
\rho (z)=\frac{\left( q^{2z},\ q^{-2z};\ q\right) _{\infty }\ \left(
q^{z}-q^{-z}\right) ^{-1}}{\left( q^{2\alpha +2z},\ q^{2\alpha -2z},\
q^{2-2\alpha +2z},\ q^{2-2\alpha -2z};\ q^{2}\right) _{\infty }},
\end{equation}
where $\alpha $ is an arbitrary additional parameter. We shall see later
that this solution satisfies the correct boundary conditions for our second
order divided-difference Askey--Wilson operator (2.10) for certain values of
this parameter $\alpha $.

Let us multiply (3.1) by $v(z)$, (3.2) by $u(z)$ and subtract the second
equality from the first one. As a result we get 
\begin{equation}
\left( \lambda -\lambda ^{\prime }\right) \ u(z)\ v(z)\ \rho (z)\nabla
x_{1}(z) 
=\Delta \left[ \sigma \rho (z)\ W\left( u(z),\ v(z)\right) \right] ,
\end{equation}
where 
\begin{eqnarray}
W\left( u(z),\ v(z)\right) &=&\left| 
\begin{array}{cc}
u(z) & v(z) \\ 
\dfrac{\nabla u(z)}{\nabla x(z)} & \dfrac{\nabla v(z)}{\nabla x(z)}
\end{array}
\right| \\
&=&u(z)\ \dfrac{\nabla v(z)}{\nabla x(z)}-v(z)\ \dfrac{\nabla u(z)}{\nabla
x(z)}  \notag
\end{eqnarray}
is the analog of the Wronskian \cite{Ni:Su:Uv}.

\begin{figure}
\unitlength=.8mm
\special{em:linewidth 0.8pt}
\linethickness{0.8pt}
\begin{picture}(158.00,158.00)
\renewcommand{\l}{\left}
\renewcommand{\r}{\right}
\renewcommand{\a}{\alpha}
\renewcommand{\b}{\beta}
\newcommand{\g}{\gamma}
\newcommand{\f}{\varphi}
\renewcommand{\t}{\theta}
\newcommand{\s}{\sigma}
\newcommand{\la}{\lambda}
\newcommand{\om}{\omega}
\newcommand{\ep}{\varepsilon}

\put(50,0){\line(0,1){158}}
\put(0,60){\line(1,0){158}}
\put(50,35){\line(1,0){70}}
\put(120,35){\line(0,1){100}}
\put(120,135){\line(-1,0){70}}
\put(70.00,60.00){\circle*{1.50}}
\put(100.00,60.00){\circle*{1.50}}
\put(70.00,115.00){\circle*{1.50}}
\put(100.00,115.00){\circle*{1.50}}
\put(150.00,65.00){\makebox(0,0)[cb]{\footnotesize Re $z$}}
\put(55.00,150.00){\makebox(0,0)[lb]{\footnotesize Im $z$}}
\put(70.00,65.00){\makebox(0,0)[cb]{\footnotesize $\alpha$}}
\put(100.00,65.00){\makebox(0,0)[cb]{\footnotesize $1-\alpha$}}
\put(70.00,110.00){\makebox(0,0)[ct]{\footnotesize $\alpha+\dfrac{i \pi}{\log q^{-1}}$}}
\put(100.00,110.00){\makebox(0,0)[ct]{\footnotesize $1-\alpha+\dfrac{i\pi}{\log q^{-1}}$}}
\put(45.00,109.00){\makebox(0,0)[rt]{\footnotesize $C$}}
\put(45.00,55.00){\makebox(0,0)[rt]{\footnotesize $0$}}
\put(125.00,55.00){\makebox(0,0)[lt]{\footnotesize $1$}}
\put(45.00,135.00){\makebox(0,0)[rc]{\footnotesize $\dfrac{3i\pi}{2\log q^{-1}}$}}
\put(45.00,35.00){\makebox(0,0)[rc]{\footnotesize $\dfrac{i\pi}{2\log q}$}}
\put(125.00,35.00){\makebox(0,0)[lc]{\footnotesize $\dfrac{i\pi}{2\log q}+1$}}
\put(125.00,135.00){\makebox(0,0)[lc]{\footnotesize $\dfrac{3i\pi}{2\log q^{-1}}+1$}}
\thicklines{
\put(50.00,35.00){\vector(0,1){100.00}}
\put(50.00,135.00){\vector(0,0){0.00}}}
\end{picture}
\end{figure}

Integrating (3.8) over the contour $C$ indicated in the Figure on the
next page; where $z$ is
such that $z=i\theta /\log q$ and $-\pi /2\leq \theta \leq 3\pi /2$; gives 
\begin{eqnarray}
&&\left( \lambda -\lambda ^{\prime }\right) \ \int_{C}u(z)\ v(z)\ \rho
(z)\nabla x_{1}(z)\ dz \\
&&\quad =\int_{C}\Delta \left[ \sigma \rho (z)\ W\left( u(z),\ v(z)\right)
\right] \ dz.  \notag
\end{eqnarray}
As a function in $z$, the integrand in the right side of (3.10) has the
natural purely imaginary period $T=2\pi i/\log q$ when $0<q<1$, so this
integral is equal to 
\begin{equation}
\int_{D}\sigma \rho (z)\ W\left( u(z),\ v(z)\right) \ dz,
\end{equation}
where $D$ is the boundary of the rectangle on the Figure oriented
counterclockwise.

The basic trigonometric functions $C(x)$ and $S(x)$ are entire functions in
the complex $z$-plane due to (2.8)--(2.9). Therefore, the poles of the
integrand in (3.11) inside the rectangle in the Figure are the simple poles
of $\rho (z)$ at $z=\alpha $, $z=1-\alpha $ and at $z=\alpha -i\pi /\log q$, 
$z=1-\alpha -i\pi /\log q$ when $0<\ $Re$\ \alpha <1/2$. Hence, by Cauchy's
theorem, 
\begin{eqnarray}
&&\frac{1}{2\pi i}\int_{D}\rho \ W\left( u,\ v\right) \ dz \\
&&\quad =\ \left. \text{Res\ }f(z)\right| _{z=\alpha }+\left. \text{Res\ }%
f(z)\right| _{z=1-\alpha }  \notag \\
&&\quad +\ \left. \text{Res\ }f(z)\right| _{z=\alpha -i\pi /\log q}+\left. 
\text{Res\ }f(z)\right| _{z=1-\alpha -i\pi /\log q},  \notag
\end{eqnarray}
where 
\begin{eqnarray}
f(z) &=&\rho (z)\ W\left( u(z),\ v(z)\right) \\
&=&\frac{q^{-z}\ \left( q^{2z},\ q^{1-2z};\ q\right) _{\infty }\ \ W\left(
u(z),\ v(z)\right) }{\left( q^{2\alpha +2z},\ q^{2\alpha -2z},\ q^{2-2\alpha
+2z},\ q^{2-2\alpha -2z};\ q^{2}\right) _{\infty }}.  \notag
\end{eqnarray}
Evaluation of the residues at these simple poles gives 
\begin{eqnarray}
\left. \text{Res\ }f(z)\right| _{z=\alpha } &=&\lim_{z\rightarrow \alpha
}\left( z-\alpha \right) \ f(z) \\
&=&-\frac{q^{-\alpha }\left( q^{2\alpha },\ q^{1-2\alpha };\ q\ \right)
_{\infty }}{2\log q^{-1}\left( q^{2},\ q^{2},\ q^{4\alpha },\ q^{2-4\alpha
};\ q^{2}\ \right) _{\infty }}  \notag \\
&&\times \left. W\left( u(z),\ v(z)\right) \right| _{z=\alpha },  \notag
\end{eqnarray}
\begin{eqnarray}
\left. \text{Res\ }f(z)\right| _{z=1-\alpha } &=&\lim_{z\rightarrow 1-\alpha
}\left( z-1+\alpha \right) \ f(z) \\
&=&-\frac{q^{-\alpha }\left( q^{2\alpha },\ q^{1-2\alpha };\ q\ \right)
_{\infty }}{2\log q^{-1}\left( q^{2},\ q^{2},\ q^{4\alpha },\ q^{2-4\alpha
};\ q^{2}\ \right) _{\infty }}  \notag \\
&&\times \left. W\left( u(z),\ v(z)\right) \right| _{z=1-\alpha },  \notag
\end{eqnarray}
\begin{eqnarray}
\left. \text{Res\ }f(z)\right| _{z=\alpha -i\pi /\log q}
&=&\lim_{z\rightarrow \alpha -i\pi /\log q}\left( z-\alpha +i\pi /\log
q\right) \ f(z) \\
&=&\frac{q^{-\alpha }\left( q^{2\alpha },\ q^{1-2\alpha };\ q\ \right)
_{\infty }}{2\log q^{-1}\left( q^{2},\ q^{2},\ q^{4\alpha },\ q^{2-4\alpha
};\ q^{2}\ \right) _{\infty }}  \notag \\
&&\times \left. W\left( u(z),\ v(z)\right) \right| _{z=\alpha -i\pi /\log q},
\notag
\end{eqnarray}
and 
\begin{eqnarray}
\left. \text{Res\ }f(z)\right| _{z=1-\alpha -i\pi /\log q}
&=&\lim_{z\rightarrow 1-\alpha }\left( z-1+\alpha +i\pi /\log q\right) \ f(z)
\\
&=&\frac{q^{-\alpha }\left( q^{2\alpha },\ q^{1-2\alpha };\ q\ \right)
_{\infty }}{2\log q^{-1}\left( q^{2},\ q^{2},\ q^{4\alpha },\ q^{2-4\alpha
};\ q^{2}\ \right) _{\infty }}  \notag \\
&&\times \left. W\left( u(z),\ v(z)\right) \right| _{z=1-\alpha -i\pi /\log
q}.  \notag
\end{eqnarray}
But 
\begin{equation}
W\left( u(z),\ v(z)\right) =\frac{v(z)u(z-1)-u(z)v(z-1)}{x(z)-x(z-1)}
\end{equation}
by (3.9) and, therefore, 
\begin{eqnarray}
&&\left. W\left( u(z),\ v(z)\right) \right| _{z=\alpha }=\left. W\left(
u(z),\ v(z)\right) \right| _{z=1-\alpha } \\
&&\ =-\left. W\left( u(z),\ v(z)\right) \right| _{z=\alpha -i\pi /\log
q}=-\left. W\left( u(z),\ v(z)\right) \right| _{z=1-\alpha -i\pi /\log q} 
\notag
\end{eqnarray}
due to the symmetries $C(x)=C(-x)$, $x(z)=x(-z)$, and $x(z)=-x(z-i\pi /\log
q)$. Thus, the residues are equal and as a result we get 
\begin{eqnarray}
&&\frac{q^{1/2}}{\left( 1-q\right) ^{2}}\left( \omega ^{2}-\omega ^{\prime
2}\right) \ \int_{C}u(z)\ v(z)\ \rho (z)\nabla x_{1}(z)\ dz \\
\qquad \quad &&\qquad \qquad \quad =-\frac{\pi i\ q^{-\alpha }\left(
q^{2\alpha },\ q^{1-2\alpha };\ q\ \right) _{\infty }}{\log q^{-1}\left(
q^{2},\ q^{2},\ q^{4\alpha },\ q^{2-4\alpha };\ q^{2}\ \right) _{\infty }} 
\notag \\
\qquad \qquad &&\qquad \qquad \qquad \quad \times \ W\left( u(\alpha ),\
v(\alpha )\right) ,  \notag
\end{eqnarray}
where $0<\ $Re$\ \alpha <1/2$.

We have established our main equation (3.20) for the case $u(z)=C\left(
x(z);\omega \right) $ and $v(z)=C\left( x(z);\omega ^{\prime }\right) $. The
same line of consideration shows that this equation is also true when $%
u(z)=S\left( x(z);\omega \right) $ and $v(z)=S\left( x(z);\omega ^{\prime
}\right) $. The corresponding analogs of the Wronskians in (3.20) can be
written as 
\begin{eqnarray}
&&W\left( C(x(z);\omega ),\ C(x(z);\omega ^{\prime })\right) \\
&&\ =\frac{2q^{1/4}}{1-q}\ \left[ \omega \ C\left( x(z);\omega ^{\prime
}\right) S\left( x(z-1/2);\omega \right) \right.  \notag \\
&&\qquad \qquad -\left. \omega ^{\prime }\ C\left( x(z);\omega \right)
S\left( x(z-1/2);\omega ^{\prime }\right) \right]  \notag
\end{eqnarray}
and 
\begin{eqnarray}
&&W\left( S(x(z);\omega ),\ S(x(z);\omega ^{\prime })\right) \\
&&\ =\frac{2q^{1/4}}{1-q}\ \left[ \omega ^{\prime }\ S\left( x(z);\omega
\right) C\left( x(z-1/2);\omega ^{\prime }\right) \right.  \notag \\
&&\qquad \qquad -\left. \omega \ S\left( x(z);\omega ^{\prime }\right)
C\left( x(z-1/2);\omega \right) \right]  \notag
\end{eqnarray}
by (2.12)--(2.13), respectively. One can see from (3.21) and (3.22) that the
right side of (3.20) vanishes in both cases when eigenvalues $\omega $ and $%
\omega ^{\prime }$ are roots of the following equation 
\begin{equation}
S_{q}\left( x(1/4);\omega \right) =S_{q}\left( x(1/4);\omega ^{\prime
}\right) =0.
\end{equation}
This is a direct analog of (1.11) for basic trigonometric functions.

In the last case, $u(z)=C\left( x(z);\omega \right) $ and $v(z)=S\left(
x(z);\omega ^{\prime }\right) $, the left side of (3.20) vanishes by
symmetry. It is interesting to verify that by using our method as well.
Equations (3.1) to (3.18) are valid again. But now 
\begin{eqnarray}
&&\left. W\left( u(z),\ v(z)\right) \right| _{z=\alpha }=\left. W\left(
u(z),\ v(z)\right) \right| _{z=1-\alpha } \\
&&\ =\left. W\left( u(z),\ v(z)\right) \right| _{z=\alpha -i\pi /\log
q}=\left. W\left( u(z),\ v(z)\right) \right| _{z=1-\alpha -i\pi /\log q} 
\notag
\end{eqnarray}
due to the symmetries $C(x)=C(-x)$, $S(x)=-S(-x)$, $x(z)=x(-z)$, and $%
x(z)=-x(z-i\pi /\log q)$. Therefore, 
\begin{eqnarray}
&&\frac{q^{1/2}}{\left( 1-q\right) ^{2}}\left( \omega ^{2}-\omega ^{\prime
2}\right) \ \int_{C}u(z)\ v(z)\ \rho (z)\nabla x_{1}(z)\ dz \\
\qquad \quad &&\qquad \qquad \quad =-\frac{\pi i\ q^{-\alpha }\left(
q^{2\alpha },\ q^{1-2\alpha };\ q\ \right) _{\infty }}{2\log q^{-1}\left(
q^{2},\ q^{2},\ q^{4\alpha },\ q^{2-4\alpha };\ q^{2}\ \right) _{\infty }} 
\notag \\
\qquad \qquad &&\qquad \qquad \qquad \quad \times \ \left[ W\left( u(\alpha
),\ v(\alpha )\right) -W\left( u(\alpha ),\ v(\alpha )\right) \right] \equiv
0,  \notag
\end{eqnarray}
when $0<\ $Re$\ \alpha <1/2$.

Combining all the above cases together, we finally arrive at the \textit{%
continuous orthogonality relations} for basic trigonometric functions, 
\begin{eqnarray}
&&\int_{0}^{\pi }C\left( \cos \theta ;\omega \right) \ C\left( \cos \theta
;\omega ^{\prime }\right) \ \frac{\ \left( e^{2i\theta },\ e^{-2i\theta };\
q\right) _{\infty }}{\left( q^{1/2}e^{2i\theta },\ q^{1/2}e^{-2i\theta };\
q\right) _{\infty }}\ d\theta \\
&&\ =\left\{ 
\begin{tabular}{lll}
$0$ & if & $\omega \neq \omega ^{\prime },$ \\ 
$\pi \dfrac{\left( q^{1/2};\ q\right) _{\infty }^{2}}{\left( q;\ q\right)
_{\infty }^{2}}\ C\left( \eta ;\omega \right) \dfrac{\partial }{\partial
\omega }S\left( \eta ;\omega \right) $ & if & $\omega =\omega ^{\prime };$%
\end{tabular}
\right.  \notag
\end{eqnarray}
\begin{eqnarray}
&&\int_{0}^{\pi }S\left( \cos \theta ;\omega \right) \ S\left( \cos \theta
;\omega ^{\prime }\right) \ \frac{\ \left( e^{2i\theta },\ e^{-2i\theta };\
q\right) _{\infty }}{\left( q^{1/2}e^{2i\theta },\ q^{1/2}e^{-2i\theta };\
q\right) _{\infty }}\ d\theta \\
&&\ =\left\{ 
\begin{tabular}{lll}
$0$ & if & $\omega \neq \omega ^{\prime },$ \\ 
$\pi \dfrac{\left( q^{1/2};\ q\right) _{\infty }^{2}}{\left( q;\ q\right)
_{\infty }^{2}}\ C\left( \eta ;\omega \right) \dfrac{\partial }{\partial
\omega }S\left( \eta ;\omega \right) $ & if & $\omega =\omega ^{\prime };$%
\end{tabular}
\right.  \notag
\end{eqnarray}
and 
\begin{equation}
\int_{0}^{\pi }C\left( \cos \theta ;\omega \right) \ S\left( \cos \theta
;\omega ^{\prime }\right) \ \frac{\ \left( e^{2i\theta },\ e^{-2i\theta };\
q\right) _{\infty }}{\left( q^{1/2}e^{2i\theta },\ q^{1/2}e^{-2i\theta };\
q\right) _{\infty }}\ d\theta =0.
\end{equation}
Here $\eta :=x(1/4)=\left( q^{1/4}+q^{-1/4}\right) /2$ and the eigenvalues $%
\omega $ and $\omega ^{\prime }$ satisfy the ``boundary'' condition (3.23).

For arbitrary $\omega \neq \omega ^{\prime }$ one gets from (3.20)--(3.22)\ 
\begin{eqnarray}
&&\int_{0}^{\pi }C\left( \cos \theta ;\omega \right) \ C\left( \cos \theta
;\omega ^{\prime }\right) \ \frac{\ \left( e^{2i\theta },\ e^{-2i\theta };\
q\right) _{\infty }}{\left( q^{1/2}e^{2i\theta },\ q^{1/2}e^{-2i\theta };\
q\right) _{\infty }}\ d\theta \\
&&\ =\frac{2\pi }{\omega ^{2}-\omega ^{\prime 2}}\ \dfrac{\left( q^{1/2};\
q\right) _{\infty }^{2}}{\left( q;\ q\right) _{\infty }^{2}}\   \notag \\
&&\ \times \left[ \omega \ C\left( \eta ;\omega ^{\prime }\right) S\left(
\eta ;\omega \right) -\omega ^{\prime }\ C\left( \eta ;\omega \right)
S\left( \eta ;\omega ^{\prime }\right) \right]  \notag
\end{eqnarray}
and 
\begin{eqnarray}
&&\int_{0}^{\pi }S\left( \cos \theta ;\omega \right) \ S\left( \cos \theta
;\omega ^{\prime }\right) \ \frac{\ \left( e^{2i\theta },\ e^{-2i\theta };\
q\right) _{\infty }}{\left( q^{1/2}e^{2i\theta },\ q^{1/2}e^{-2i\theta };\
q\right) _{\infty }}\ d\theta \\
&&\ =\frac{2\pi }{\omega ^{2}-\omega ^{\prime 2}}\ \dfrac{\left( q^{1/2};\
q\right) _{\infty }^{2}}{\left( q;\ q\right) _{\infty }^{2}}\   \notag \\
&&\ \times \left[ \omega ^{\prime }\ S\left( \eta ;\omega \right) C\left(
\eta ;\omega ^{\prime }\right) -\omega \ S\left( \eta ;\omega ^{\prime
}\right) C\left( \eta ;\omega \right) \right] .  \notag
\end{eqnarray}
Also, in the limit $\omega \rightarrow \omega ^{\prime }$, 
\begin{eqnarray}
&&\int_{0}^{\pi }C^{2}\left( \cos \theta ;\omega \right) \ \frac{\ \left(
e^{2i\theta },\ e^{-2i\theta };\ q\right) _{\infty }}{\left(
q^{1/2}e^{2i\theta },\ q^{1/2}e^{-2i\theta };\ q\right) _{\infty }}\ d\theta
\\
&&\ =\dfrac{\pi \left( q^{1/2};\ q\right) _{\infty }^{2}}{\omega \left( q;\
q\right) _{\infty }^{2}}\ \ \left[ \omega \ C\left( \eta ;\omega \right) 
\frac{\partial }{\partial \omega }S\left( \eta ;\omega \right) \right. 
\notag \\
&&\ \left. \quad +\ C\left( \eta ;\omega \right) S\left( \eta ;\omega
\right) -\omega \frac{\partial }{\partial \omega }C\left( \eta ;\omega
\right) S\left( \eta ;\omega \right) \right]  \notag
\end{eqnarray}
and 
\begin{eqnarray}
&&\int_{0}^{\pi }S^{2}\left( \cos \theta ;\omega \right) \ \frac{\ \left(
e^{2i\theta },\ e^{-2i\theta };\ q\right) _{\infty }}{\left(
q^{1/2}e^{2i\theta },\ q^{1/2}e^{-2i\theta };\ q\right) _{\infty }}\ d\theta
\\
&&\ =\dfrac{\pi \left( q^{1/2};\ q\right) _{\infty }^{2}}{\omega \left( q;\
q\right) _{\infty }^{2}}\ \ \left[ \omega \ C\left( \eta ;\omega \right) 
\frac{\partial }{\partial \omega }S\left( \eta ;\omega \right) \right. 
\notag \\
&&\ \left. \quad -\ C\left( \eta ;\omega \right) S\left( \eta ;\omega
\right) -\omega \frac{\partial }{\partial \omega }C\left( \eta ;\omega
\right) S\left( \eta ;\omega \right) \right] .  \notag
\end{eqnarray}
We remind the reader that $\eta $ is defined by $\eta =x(1/4)=\left(
q^{1/4}+q^{-1/4}\right) /2$. This notation will be used throughout this work.

\section{Formal Limit $q\rightarrow 1^{-}$}

In this section we formally obtain orthogonality of the trigonometric
functions as limiting cases of our orthogonality relations (3.26)--(3.28)
for basic trigonometric functions. According to \cite{Su}, 
\begin{eqnarray}
\lim_{q\rightarrow 1^{-}}C_{q}\left( x,y;\ \frac{1}{2}\omega \left(
1-q\right) \right) &=&\cos \omega \left( x+y\right) , \\
\lim_{q\rightarrow 1^{-}}S_{q}\left( x,y;\ \frac{1}{2}\omega \left(
1-q\right) \right) &=&\sin \omega \left( x+y\right) .
\end{eqnarray}
If $\omega \neq \omega ^{\prime }$ we can rewrite (3.26) as 
\begin{equation}
\int_{0}^{\pi }C\left( \cos \theta ;\omega \right) \ C\left( \cos \theta
;\omega ^{\prime }\right) \ \left( e^{2i\theta },e^{-2i\theta };q\right)
_{1/2}\ d\theta =0,
\end{equation}
where 
\begin{equation}
\left( a;r\right) _{\alpha }:=\frac{\left( a;r\right) _{\infty }}{\left(
ar^{\alpha };r\right) _{\infty }}.
\end{equation}
Using the limiting relation \cite{Ga:Ra} 
\begin{equation}
\lim_{q\rightarrow 1^{-}}\left( a;r\right) _{\alpha }=\left( 1-a\right)
^{\alpha },
\end{equation}
one can see that 
\begin{equation}
\left( e^{2i\theta },e^{-2i\theta };q\right) _{1/2}\rightarrow 2\sin \theta
\end{equation}
as $q\rightarrow 1^{-}$. Therefore, changing $\omega $ to $\left( 1-q\right)
\omega /2$ in (4.3), with the help of (4.1) when $y=0$ we obtain the
orthogonality relation (1.6) with $l=1$. The boundary condition (1.11)
follows from (3.23) in the same limit.

When $\omega =\omega ^{\prime }$ we can rewrite (3.26) as 
\begin{eqnarray}
\int_{0}^{\pi }C^{2}\left( \cos \theta ;\omega \right) \ \left( e^{2i\theta
},e^{-2i\theta };q^{2}\right) _{1/2}\ d\theta \ && \\
\ =\frac{\pi \left( 1-q\right) }{\Gamma _{q}^{2}\left( 1/2\right) }\ C\left(
\eta ;\omega \right) \frac{\partial }{\partial \omega }S\left( \eta ;\omega
\right) ,\ &&  \notag
\end{eqnarray}
where $\eta =\left( q^{1/4}+q^{-1/4}\right) /2$ and 
\begin{equation}
\Gamma _{q}\left( z\right) =\left( 1-q\right) ^{1-z}\ \frac{\left(
q;q\right) _{\infty }}{\left( q^{z};q\right) _{\infty }}
\end{equation}
is a $q$-analog of Euler's gamma function $\Gamma \left( z\right) $ (see,
for example, \cite{Ga:Ra}). Changing $\omega $ to $\left( 1-q\right) \omega
/2$ in (4.8), with the aid of 
\begin{equation}
\lim_{q\rightarrow 1^{-}}\Gamma _{q}\left( z\right) =\Gamma \left( z\right) ,
\end{equation}
we get 
\begin{equation}
2\int_{-1}^{1}\cos ^{2}\pi nx\ dx=\frac{2\pi }{\Gamma ^{2}\left( 1/2\right) }%
\ \cos ^{2}\pi n=2,
\end{equation}
where $n=\pm 1,\pm 2,...$ , in the limit $q\rightarrow 1^{-}$.

In a similar manner one can obtain (1.7) and (1.8) from (3.27) and (3.28),
respectively.

\section{Some Properties of Zeros}

In Section 3 we have established the orthogonality relations for the basic
trigonometric functions (3.26)--(3.28) under the boundary condition (3.23).
Here we would like to discuss some properties of $\omega $-zeros of the
corresponding basic sine function, 
\begin{eqnarray}
S\left( \eta ;\omega \right) &=&\frac{\left( -\omega ^{2};q^{2}\right)
_{\infty }}{\left( -q\omega ^{2};q^{2}\right) _{\infty }}\ \frac{\omega }{%
1-q^{1/2}}\  \\
&&\times \ _{2}\varphi _{1}\left( 
\begin{array}{c}
\begin{array}{ll}
-q^{3/2}, & -q^{5/2}
\end{array}
\\ 
q^{3}
\end{array}
;\ q^{2},\ -\omega ^{2}\right)  \notag \\
&=&\frac{\left( q^{3/2}\omega ^{2};\ q\right) _{\infty }}{\left(
q^{3},-q\omega ^{2};\ q^{2}\right) _{\infty }}\ \frac{\omega }{1-q^{1/2}} 
\notag \\
&&\times \ _{2}\varphi _{2}\left( 
\begin{array}{cc}
-\omega ^{2}, & -q\omega ^{2} \\ 
q^{3/2}\omega ^{2}, & q^{5/2}\omega ^{2}
\end{array}
;\ q^{2},\ q^{3}\right) ,  \notag
\end{eqnarray}
and the basic cosine function, 
\begin{eqnarray}
C\left( \eta ;\omega \right) &=&\frac{\left( -\omega ^{2};q^{2}\right)
_{\infty }}{\left( -q\omega ^{2};q^{2}\right) _{\infty }}\  \\
&&\times \ _{2}\varphi _{1}\left( 
\begin{array}{c}
\begin{array}{ll}
-q^{1/2}, & -q^{3/2}
\end{array}
\\ 
q
\end{array}
;\ q^{2},\ -\omega ^{2}\right)  \notag \\
&=&\frac{\left( q^{1/2}\omega ^{2};\ q\right) _{\infty }}{\left( q,-q\omega
^{2};\ q^{2}\right) _{\infty }}  \notag \\
&&\times \ _{2}\varphi _{2}\left( 
\begin{array}{cc}
-\omega ^{2}, & -q\omega ^{2} \\ 
q^{1/2}\omega ^{2}, & q^{3/2}\omega ^{2}
\end{array}
;\ q^{2},\ q\right) .  \notag
\end{eqnarray}
One can see that these functions have almost the same structure as the $q$%
-Bessel function discussed in \cite{Is:Ma:Su1}, \cite{Is:Ma:Su2}. So we can
apply a similar method to establish main properties of zeros of the
functions (5.1)--(5.2).

The first property is that the $q$-sine function $S\left( \eta ;\omega
\right) $ has an infinity of real $\omega $-zeros. To prove that we can
again consider the large $\omega $-asymptotics of the function (5.1). The $%
_{2}\varphi _{1}$ here can be transformed by (III.1) of \cite{Ga:Ra}, which
gives 
\begin{eqnarray}
S\left( \eta ;\omega \right) &=&\frac{\left( -q^{5/2},q^{3/2}\omega
^{2};q^{2}\right) _{\infty }}{\left( q^{3},-q\omega ^{2};q^{2}\right)
_{\infty }}\ \frac{\omega }{1-q^{1/2}} \\
&&\times \ _{2}\varphi _{1}\left( 
\begin{array}{c}
\begin{array}{ll}
-q^{1/2}, & -\omega ^{2}
\end{array}
\\ 
q^{3/2}\omega ^{2}
\end{array}
;\ q^{2},\ -q^{5/2}\right) .  \notag
\end{eqnarray}
For large values of $\omega $, such that $\omega ^{2}\neq q^{-3/2-2n}$ where 
$n=0,1,2,...\ $, 
\begin{eqnarray}
&&_{2}\varphi _{1}\left( 
\begin{array}{c}
\begin{array}{ll}
-q^{1/2}, & -\omega ^{2}
\end{array}
\\ 
q^{3/2}\omega ^{2}
\end{array}
;\ q^{2},\ -q^{5/2}\right) \\
&&\ \rightarrow \ _{1}\varphi _{0}\left( 
\begin{array}{c}
-q^{1/2} \\ 
-
\end{array}
;\ q^{2},\ q\right) =\frac{\left( -q^{3/2};q^{2}\right) _{\infty }}{\left(
q;q^{2}\right) _{\infty }},  \notag
\end{eqnarray}
by the $q$-binomial theorem. Therefore, as $\omega \rightarrow \infty $, 
\begin{eqnarray}
S\left( \eta ;\omega \right) &=&\frac{\left( -q^{1/2};q\right) _{\infty }}{%
\left( q;q^{2}\right) _{\infty }^{2}} \\
&&\times \ \omega \ \frac{\left( q^{3/2}\omega ^{2};q^{2}\right) _{\infty }}{%
\left( -q\omega ^{2};q^{2}\right) _{\infty }}\left[ 1+\text{o}(1)\right] , 
\notag
\end{eqnarray}
by (5.3) and (5.4). But the function 
\begin{equation*}
\left( q^{3/2}\omega ^{2};q^{2}\right) _{\infty }
\end{equation*}
oscillates and has an infinity of real zeros as $\omega $ approaches
infinity. Indeed, consider the points $\omega =\gamma _{n}$, such that 
\begin{equation}
\gamma _{n}^{2}=\beta ^{2}q^{-2n},
\end{equation}
where $n=0,1,2,...$ and $q^{1/2}<\beta ^{2}<q^{-3/2}$, as test points. Then,
by using (I.9) of \cite{Ga:Ra}, 
\begin{eqnarray}
S\left( \eta ;\gamma _{n}\right) &=&\frac{\left( -q^{1/2};q\right) _{\infty }%
}{\left( q;q^{2}\right) _{\infty }^{2}}\ \beta \ \frac{\left( q^{3/2}\beta
^{2};q^{2}\right) _{\infty }}{\left( -q\beta ^{2};q^{2}\right) _{\infty }} \\
&&\times \left( -1\right) ^{n}\ q^{-n/2}\ \frac{\left( q^{1/2}/\beta
^{2};q^{2}\right) _{n}}{\left( -q/\beta ^{2};q^{2}\right) _{n}}\left[ 1+%
\text{o}(1)\right] ,  \notag
\end{eqnarray}
as $n\rightarrow \infty $, and one can see that the right side of (5.7)
changes sign infinitely many times at the test points $\omega =\gamma _{n}$
as $\omega $ approaches infinity.

In a similar manner, one can prove that the $q$-cosine function $C\left(
\eta ;\omega \right) $ has an infinity of real $\omega $-zeros also.

Thus we have established the following theorem.

\begin{theorem}
The basic sine $S\left( \eta ;\omega \right) $ and basic cosine $C\left(
\eta ;\omega \right) $ functions have an infinity of real $\omega $-zeros
when $0<q<1$.
\end{theorem}

Now we can prove our next result.

\begin{theorem}
The basic sine $S\left( \eta ;\omega \right) $ and basic cosine $C\left(
\eta ;\omega \right) $ functions have only real $\omega $-zeros when $0<q<1$.
\end{theorem}

\proof%
Suppose that $\omega _{0}$ is a zero of the basic sine function (5.1) which
is not real. It follows from (5.1) and (III.4) of \cite{Ga:Ra} that 
\begin{eqnarray}
S\left( \eta ;\omega \right) &=&\frac{\left( q^{5/2}\omega ^{2};q^{2}\right)
_{\infty }}{\left( -q\omega ^{2};q^{2}\right) _{\infty }}\ \frac{\omega }{%
1-q^{1/2}}\  \\
&&\times \ _{2}\varphi _{2}\left( 
\begin{array}{cc}
-q^{3/2}, & -q^{5/2} \\ 
q^{3}, & q^{5/2}\omega ^{2}
\end{array}
;\ q^{2},\ q^{3/2}\omega ^{2}\right) .  \notag
\end{eqnarray}
Now we can see that $\omega _{0}$ is not purely imaginary, because otherwise
our function would be a multiple of a positive function.

Let $\omega _{1}$ be the complex number conjugate to $\omega _{0}$, so that $%
\omega _{1}$ is also a zero of (5.1) because this function is a real
function of $\omega $. Since $\omega _{0}^{2}\neq $ $\omega _{1}^{2}$ the
integral in the orthogonality relation (3.26) equals zero, but the integrand
on the left is positive, and so we have obtained a contradiction. Hence a
complex zero $\omega _{0}$ cannot exist. One can consider the case of the
basic cosine function in a similar fashion. 
\endproof%

\begin{theorem}
If $0<q<1$, then the real $\omega $-zeros of the basic sine $S\left( \eta
;\omega \right) $ and basic cosine $C\left( \eta ;\omega \right) $ functions
are simple.
\end{theorem}

\proof%
This follows directly from the relations (3.31) and (3.32). Consider, for
example, the case of the basic sine function. If $\omega =\omega ^{\prime }$%
, then the integral in the left side of (3.31) is positive, which means that 
$\dfrac{\partial }{\partial \omega }S\left( \eta ;\omega \right) \neq 0$
when $S\left( \eta ;\omega \right) =0$. The same is true for the zeros of
the basic cosine function. 
\endproof%

Our next property is that the positive zeros of the basic sine function $%
S\left( \eta ;\omega \right) $ are interlaced with those of the basic cosine
function $C\left( \eta ;\omega \right) $.

\begin{theorem}
If $\omega _{1},\omega _{2},\omega _{3},...$ are the positive zeros of $%
S\left( \eta ;\omega \right) $ arranged in ascending order of magnitude, and 
$\varpi _{1},\varpi _{2},\varpi _{3},...$ are those of $C\left( \eta ;\omega
\right) $, then 
\begin{equation}
0=\omega _{0}<\varpi _{1}<\omega _{1}<\varpi _{2}<\omega _{3}<\varpi
_{3}<...\ ,
\end{equation}
if $0<q<1$.
\end{theorem}

\proof%
Suppose that $\omega _{k}$ and $\omega _{k+1}$ are two successive zeros of $%
S\left( \eta ;\omega \right) $. Then the derivative $\dfrac{\partial }{%
\partial \omega }S\left( \eta ;\omega \right) $ has different signs at $%
\omega =\omega _{k}$ and $\omega =\omega _{k+1}$. This means, in view of
(3.32), that $C\left( \eta ;\omega \right) $ changes its sign between $%
\omega _{k}$ and $\omega _{k+1}$ and , therefore, has at least one zero on
each interval $\left( \omega _{k},\ \omega _{k+1}\right) $.

To complete the proof of the theorem, we should show that $C\left( \eta
;\omega \right) $ changes its sign on each interval $\left( \omega _{k},\
\omega _{k+1}\right) $ only once. Suppose that $C\left( \eta ;\varpi
_{k}\right) =$ $C\left( \eta ;\varpi _{k+1}\right) =0$ and $\omega
_{k}<\varpi _{k}<\varpi _{k+1}<\omega _{k+1}$. Then, by (3.32), the function 
$S\left( \eta ;\omega \right) $ has different signs at $\omega =\varpi _{k}$
and $\omega =\varpi _{k+1}$ and, therefore, this function has at least one
more zero on $\left( \omega _{k},\ \omega _{k+1}\right) $. So, we have
obtained a contradiction, and, therefore, the basic cosine function $C\left(
\eta ;\omega \right) $ has exactly one zero between any two successive zeros
of the basic sine function $S\left( \eta ;\omega \right) $. 
\endproof%

The proof of Theorem 5.1 has strongly indicated that asymptotically the
large $\omega $-zeros of the basic sine function $S\left( \eta ;\omega
\right) $ are 
\begin{equation}
\omega _{n}=\pm \varkappa _{n}\ q^{-n},\qquad q^{1/4}\leq \varkappa
_{n}<q^{-3/4}
\end{equation}
as $n\rightarrow \infty $. The same consideration as in \cite{Is} and \cite
{Is:Ma:Su2} shows that $S\left( \eta ;\omega \right) $ changes sign only
once between any two successive test points $\omega =\gamma _{n}$ and $%
\omega =\gamma _{n+1}$ determined by (5.6) for large values of $n$. We
include details of this proof in Section 16 to make this work as
self-contained as possible.

Our next theorem provides a more accurate estimate for the distribution of
the large zeros of this function.

\begin{theorem}
If $\omega _{1},\omega _{2},\omega _{3},...$ are the positive zeros of $%
S\left( \eta ;\omega \right) $ arranged in ascending order of magnitude,
then 
\begin{equation}
\omega _{n}=q^{1/4-n}+\text{o}\left( 1\right) ,
\end{equation}
as $n\rightarrow \infty $.
\end{theorem}

\proof%
In view of (5.1) and (III.32) of \cite{Ga:Ra}, 
\begin{eqnarray}
S\left( \eta ;\omega \right) &=&\ \frac{\omega }{1-q^{1/2}}\ \frac{\left(
-q^{3/2},-q^{5/2},q^{3/2}\omega ^{2},q^{1/2}/\omega ^{2};q^{2}\right)
_{\infty }}{\left( q,q^{3},-q\omega ^{2},-q^{2}/\omega ^{2};q^{2}\right)
_{\infty }} \\
&&\quad \quad \times \ _{2}\varphi _{1}\left( 
\begin{array}{c}
\begin{array}{ll}
-q^{1/2}, & -q^{3/2}
\end{array}
\\ 
q
\end{array}
;\ q^{2},\ -\frac{q}{\omega ^{2}}\right)  \notag \\
&&+\ \frac{\omega }{1-q^{1/2}}\ \frac{\left( -q^{1/2},-q^{3/2},q^{5/2}\omega
^{2},q^{-1/2}/\omega ^{2};q^{2}\right) _{\infty }}{\left(
q^{-1},q^{3},-q\omega ^{2},-q^{2}/\omega ^{2};q^{2}\right) _{\infty }} 
\notag \\
&&\quad \quad \quad \times \ _{2}\varphi _{1}\left( 
\begin{array}{c}
\begin{array}{ll}
-q^{3/2}, & -q^{5/2}
\end{array}
\\ 
q^{3}
\end{array}
;\ q^{2},\ -\frac{q}{\omega ^{2}}\right) ,  \notag
\end{eqnarray}
which gives the large $\omega $-asymptotic of $S\left( \eta ;\omega \right) $%
. When $\omega =q^{1/4-n}$ and $n=1,2,3,...,$ the first term in (5.12)
vanishes and we get 
\begin{eqnarray}
S\left( \eta ;q^{1/4-n}\right) &=&\left( -1\right) ^{n}q^{n/2-1/4}\ \frac{%
1+q^{1/2}}{1-q^{1/2}}\ \frac{\left( -q^{3/2};q^{2}\right) _{n}}{\left(
-q^{1/2};q^{2}\right) _{n}} \\
&&\times \ _{2}\varphi _{1}\left( 
\begin{array}{c}
\begin{array}{ll}
-q^{3/2}, & -q^{5/2}
\end{array}
\\ 
q^{3}
\end{array}
;\ q^{2},\ -q^{2n+1/2}\right)  \notag
\end{eqnarray}
with the help of (I.9) of \cite{Ga:Ra}, Thus, 
\begin{equation}
\lim_{n\rightarrow \infty }S\left( \eta ;q^{1/4-n}\right) =0,
\end{equation}
which proves our theorem. 
\endproof%

In a similar fashion, one can establish the following theorem.

\begin{theorem}
If $\varpi _{1},\varpi _{2},\varpi _{3},...$ are the positive zeros of $%
C\left( \eta ;\omega \right) $ arranged in ascending order of magnitude,
then 
\begin{equation}
\varpi _{n}=q^{3/4-n}+\text{o}\left( 1\right) ,
\end{equation}
as $n\rightarrow \infty $.
\end{theorem}

The asymptotic formulas (5.11) and (5.15) for large $\omega $-zeros of the
basic sine $S\left( \eta ;\omega \right) $ and basic cosine $C\left( \eta
;\omega \right) $ functions confirm the interlacing property (5.9) from
Theorem 5.4.

Let us also discuss the large $\omega $-asymptotics of the basic sine $%
S\left( x;\omega \right) $ and basic cosine $C\left( x;\omega \right) $
functions when $x=\cos \theta $ belongs to the interval of orthogonality $%
-1<x<1$. From (2.1) and (2.2) one gets 
\begin{eqnarray}
C\left( \cos \theta ;\omega \right) &=&\ _{2}\varphi _{1}\left( 
\begin{array}{c}
\begin{array}{ll}
-e^{2i\theta }, & -e^{-2i\theta }
\end{array}
\\ 
q
\end{array}
;\ q^{2},\ -q\omega ^{2}\right) \\
&=&\frac{\left( -e^{-2i\theta },-qe^{-2i\theta },\ q\omega ^{2}e^{2i\theta
},\ e^{-2i\theta }/\omega ^{2};q^{2}\right) _{\infty }}{\left( q,\
e^{-4i\theta },-q\omega ^{2},-q/\omega ^{2};q^{2}\right) _{\infty }}  \notag
\\
&&\times \ _{2}\varphi _{1}\left( 
\begin{array}{c}
\begin{array}{ll}
-e^{2i\theta }, & -qe^{2i\theta }
\end{array}
\\ 
q^{2}e^{4i\theta }
\end{array}
;\ q^{2},\ -\frac{q^{2}}{\omega ^{2}}\right)  \notag \\
&&+\frac{\left( -e^{2i\theta },-qe^{2i\theta },\ q\omega ^{2}e^{-2i\theta
},\ e^{2i\theta }/\omega ^{2};q^{2}\right) _{\infty }}{\left( q,\
e^{4i\theta },-q\omega ^{2},-q/\omega ^{2};q^{2}\right) _{\infty }}  \notag
\\
&&\times \ _{2}\varphi _{1}\left( 
\begin{array}{c}
\begin{array}{ll}
-e^{-2i\theta }, & -qe^{-2i\theta }
\end{array}
\\ 
q^{2}e^{-4i\theta }
\end{array}
;\ q^{2},\ -\frac{q^{2}}{\omega ^{2}}\right)  \notag
\end{eqnarray}
and 
\begin{eqnarray}
S\left( \cos \theta ;\omega \right) &=&\ \frac{2q^{1/4}\omega }{1-q}\ \cos
\theta \\
&&\times \ _{2}\varphi _{1}\left( 
\begin{array}{c}
\begin{array}{ll}
-qe^{2i\theta }, & -qe^{-2i\theta }
\end{array}
\\ 
q^{3}
\end{array}
;\ q^{2},\ -q\omega ^{2}\right)  \notag \\
&=&\frac{2q^{1/4}\omega }{1-q}\ \cos \theta  \notag \\
&&\times \Biggl[\frac{\left( -qe^{-2i\theta },-q^{2}e^{-2i\theta },\
q^{2}\omega ^{2}e^{2i\theta },\ e^{-2i\theta }/\omega ^{2};q^{2}\right)
_{\infty }}{\left( q^{3},\ e^{-4i\theta },-q\omega ^{2},-q/\omega
^{2};q^{2}\right) _{\infty }}  \notag \\
&&\quad \times \ _{2}\varphi _{1}\left( 
\begin{array}{c}
\begin{array}{ll}
-e^{2i\theta }, & -qe^{2i\theta }
\end{array}
\\ 
q^{2}e^{4i\theta }
\end{array}
;\ q^{2},\ -\frac{q^{2}}{\omega ^{2}}\right)  \notag \\
&&\quad +\frac{\left( -qe^{2i\theta },-q^{2}e^{2i\theta },\ q^{2}\omega
^{2}e^{-2i\theta },\ e^{2i\theta }/\omega ^{2};q^{2}\right) _{\infty }}{%
\left( q^{3},\ e^{4i\theta },-q\omega ^{2},-q/\omega ^{2};q^{2}\right)
_{\infty }}  \notag \\
&&\quad \times \ _{2}\varphi _{1}\left( 
\begin{array}{c}
\begin{array}{ll}
-e^{-2i\theta }, & -qe^{-2i\theta }
\end{array}
\\ 
q^{2}e^{-4i\theta }
\end{array}
;\ q^{2},\ -\frac{q^{2}}{\omega ^{2}}\right) \Biggr]  \notag
\end{eqnarray}
by (III.3) and (III.32) of \cite{Ga:Ra}. For $|x|<1$, $|q|<1$ and large $%
\omega $ it is clear from (5.16) and (5.17) that the leading terms in the
asymptotic expansions of $C\left( \cos \theta ;\omega \right) $ and $S\left(
\cos \theta ;\omega \right) $ are given by 
\begin{eqnarray}
C\left( \cos \theta ;\omega \right) &\sim &\frac{\left( -e^{-2i\theta
};q\right) _{\infty }}{\left( q,\ e^{-4i\theta };q^{2}\right) _{\infty }}\ 
\frac{\left( q\omega ^{2}e^{2i\theta };q^{2}\right) _{\infty }}{\left(
-q\omega ^{2};q^{2}\right) _{\infty }} \\
&&+\ \frac{\left( -e^{2i\theta };q\right) _{\infty }}{\left( q,\ e^{4i\theta
};q^{2}\right) _{\infty }}\ \frac{\left( q\omega ^{2}e^{-2i\theta
};q^{2}\right) _{\infty }}{\left( -q\omega ^{2};q^{2}\right) _{\infty }} 
\notag
\end{eqnarray}
and 
\begin{eqnarray}
S\left( \cos \theta ;\omega \right) &\sim &\frac{2q^{1/4}\omega }{1-q}\ \cos
\theta \\
&&\times \left[ \frac{\left( -qe^{-2i\theta };q\right) _{\infty }}{\left(
q^{3},\ e^{-4i\theta };q^{2}\right) _{\infty }}\ \frac{\left( q^{2}\omega
^{2}e^{2i\theta };q^{2}\right) _{\infty }}{\left( -q\omega ^{2};q^{2}\right)
_{\infty }}\right.  \notag \\
&&\quad +\left. \frac{\left( -qe^{2i\theta };q\right) _{\infty }}{\left(
q^{3},\ e^{4i\theta };q^{2}\right) _{\infty }}\ \frac{\left( q^{2}\omega
^{2}e^{-2i\theta };q^{2}\right) _{\infty }}{\left( -q\omega
^{2};q^{2}\right) _{\infty }}\right] ,  \notag
\end{eqnarray}
respectively. In particular, when $\omega =\omega _{n}$ are large zeros of
the basic sine function $S\left( \eta ;\omega \right) $ we can estimate 
\begin{equation}
C\left( \cos \theta ;\omega _{n}\right) \sim C\left( \cos \theta
;q^{1/4-n}\right) ,
\end{equation}
\begin{equation}
S\left( \cos \theta ;\omega _{n}\right) \sim S\left( \cos \theta
;q^{1/4-n}\right)
\end{equation}
due to (5.11) as $n\rightarrow \infty $. Relations (5.18)--(5.21) lead to
the following theorem.

\begin{theorem}
For $-1<x=\cos \theta <1$ and $|q|<1$ the leading term in the asymptotic
expansion of $C\left( \cos \theta ;\omega _{n}\right) $ as $n\rightarrow
\infty $ is given by 
\begin{eqnarray}
C\left( \cos \theta ;q^{1/4-n}\right) &\sim &2\frac{\left( q^{1/2};q\right)
_{\infty }}{\left( q;q^{2}\right) _{\infty }^{2}} \\
&&\times \left| A\left( e^{i\theta }\right) \right| \ \cos \left( \left(
2\theta +\pi \right) n-\chi \right) ,  \notag
\end{eqnarray}
where 
\begin{equation}
A\left( e^{i\theta }\right) =\left( 1-q^{1/2}e^{2i\theta }\right) \ \frac{%
\left( q^{3/2}e^{-2i\theta },q^{5/2}e^{2i\theta };q^{2}\right) _{\infty }}{%
\left( e^{2i\theta };q\right) _{\infty }},
\end{equation}
\begin{equation}
\left| A\left( e^{i\theta }\right) \right| ^{-2}=\frac{\left( e^{2i\theta
},\ e^{-2i\theta };q\right) _{\infty }}{\left( q^{1/2}e^{2i\theta },\
q^{1/2}e^{-2i\theta };q\right) _{\infty }},
\end{equation}
and 
\begin{equation}
\chi =\arg A\left( e^{i\theta }\right) .
\end{equation}

For $-1<x=\cos \theta <1$ and $|q|<1$ the leading term in the asymptotic
expansion of $S\left( \cos \theta ;\omega _{n}\right) $ as $n\rightarrow
\infty $ is given by 
\begin{eqnarray}
S\left( \cos \theta ;q^{1/4-n}\right) &\sim &2\frac{\left( q^{1/2};q\right)
_{\infty }}{\left( q;q^{2}\right) _{\infty }^{2}} \\
&&\times \left| B\left( e^{i\theta }\right) \right| \ \cos \left( \left(
2\theta +\pi \right) \left( n-1\right) -\psi \right) ,  \notag
\end{eqnarray}
where 
\begin{equation}
B\left( e^{i\theta }\right) =e^{i\theta }\ \frac{\left( q^{1/2}e^{-2i\theta
},q^{3/2}e^{2i\theta };q^{2}\right) _{\infty }}{\left( e^{2i\theta
};q\right) _{\infty }},
\end{equation}
\begin{equation}
\left| B\left( e^{i\theta }\right) \right| ^{-2}=\frac{\left( e^{2i\theta
},\ e^{-2i\theta };q\right) _{\infty }}{\left( q^{1/2}e^{2i\theta },\
q^{1/2}e^{-2i\theta };q\right) _{\infty }},
\end{equation}
and 
\begin{equation}
\psi =\arg B\left( e^{i\theta }\right) .
\end{equation}
\end{theorem}

From (5.23) and (5.27), 
\begin{equation}
A\left( e^{i\theta }\right) =e^{i\theta }\ \frac{\left( e^{-2i\theta
};q\right) _{\infty }}{\left( e^{2i\theta };q\right) _{\infty }}\ B\left(
e^{-i\theta }\right) .
\end{equation}
It is worth mentioning also that the factor $\left| A\left( e^{i\theta
}\right) \right| ^{-2}=\left| B\left( e^{i\theta }\right) \right| ^{-2}$
coinsides with the weight function in our orthogonality relations
(3.26)--(3.28) for the basic trigonometric functions.

In a similar fashion, one can use the first lines in (5.16), (5.17), and
Exercise 3.8 of \cite{Ga:Ra} (see also the same line of reasonings in \cite
{Ch:Is:Mu}) to establish complete asymptotic expansions of the basic sine
and cosine functions for the large values of $\omega $.

\begin{theorem}
For $-1<x=\cos \theta <1$ and $|q|<1$ complete asymptotic expansions of $%
C\left( \cos \theta ;\omega \right) $ and $S\left( \cos \theta ;\omega
\right) $ as $\left| \omega \right| \rightarrow \infty $ are given by 
\begin{eqnarray}
C\left( \cos \theta ;\omega \right) &=&\frac{\left( q\omega ^{2}e^{2i\theta
};q^{2}\right) _{\infty }}{\left( e^{-2i\theta };q\right) _{\infty }\left(
q,-q\omega ^{2};q^{2}\right) _{\infty }} \\
&&\times \sum_{n=0}^{\infty }q^{2n}\frac{\left( -e^{2i\theta };q\right) _{2n}%
}{\left( q^{2},q^{2}e^{4i\theta };q^{2}\right) _{n}}\ \left( q\omega
^{2}e^{2i\theta };q^{2}\right) _{n}^{-1}  \notag \\
&&+\frac{\left( q\omega ^{2}e^{-2i\theta };q^{2}\right) _{\infty }}{\left(
e^{2i\theta };q\right) _{\infty }\left( q,-q\omega ^{2};q^{2}\right)
_{\infty }}  \notag \\
&&\times \sum_{n=0}^{\infty }q^{2n}\frac{\left( -e^{-2i\theta };q\right)
_{2n}}{\left( q^{2},q^{2}e^{-4i\theta };q^{2}\right) _{n}}\ \left( q\omega
^{2}e^{-2i\theta };q^{2}\right) _{n}^{-1}  \notag
\end{eqnarray}
and 
\begin{eqnarray}
S\left( \cos \theta ;\omega \right) &=&e^{i\theta }\frac{\left( q^{2}\omega
^{2}e^{2i\theta };q^{2}\right) _{\infty }}{\left( e^{-2i\theta };q\right)
_{\infty }\left( q,-q\omega ^{2};q^{2}\right) _{\infty }} \\
&&\times \sum_{n=0}^{\infty }q^{2n+1/4}\frac{\left( -qe^{2i\theta };q\right)
_{2n}}{\left( q^{2},q^{2}e^{4i\theta };q^{2}\right) _{n}}\ \left(
q^{2}\omega ^{2}e^{2i\theta };q^{2}\right) _{n}^{-1}  \notag \\
&&+e^{-i\theta }\frac{\left( q^{2}\omega ^{2}e^{-2i\theta };q^{2}\right)
_{\infty }}{\left( e^{2i\theta };q\right) _{\infty }\left( q,-q\omega
^{2};q^{2}\right) _{\infty }}  \notag \\
&&\times \sum_{n=0}^{\infty }q^{2n+1/4}\frac{\left( -qe^{-2i\theta
};q\right) _{2n}}{\left( q^{2},q^{2}e^{-4i\theta };q^{2}\right) _{n}}\
\left( q\omega ^{2}e^{-2i\theta };q^{2}\right) _{n}^{-1}.  \notag
\end{eqnarray}
\end{theorem}

The asymptotic expansions (5.31)--(5.32) are not in terms of the usual
asymptotic sequence $\left\{ \left( x\omega \right) ^{-n}\right\}
_{n=0}^{\infty }$, but are sums of two complete asymptotic expansions in
terms of the ``inverse generalized powers'' $\left( q^{2}\omega ^{2}e^{\pm
2i\theta };q^{2}\right) _{n}^{-1}$(cf. \cite{Ch:Is:Mu}).

\begin{remark}
{\rm Mourad Ismail pointed out to our attention the following quadratic
transformation formula 
\begin{eqnarray}
&&_{2}\varphi _{1}\left( 
\begin{array}{c}
\begin{array}{ll}
-q^{\nu +1}, & -q^{\nu +2}
\end{array}
\\ 
q^{2\nu +2}
\end{array}
;\ q^{2},\ -r^{2}/4\right) \\
&&\ =\frac{\left( q;q\right) _{\infty }}{\left( q^{\nu +1};q\right) _{\infty
}}\ \frac{\left( 2/r\right) ^{\nu }}{\left( -r^{2}/4;q^{2}\right) _{\infty }}%
\ J_{\nu }^{\left( 2\right) }\left( r;q\right) ,  \notag
\end{eqnarray}
where $|r|<2$, relating the $_{2}\varphi _{1}$ of a given structure with
Jackson's basic Bessel functions $J_{\nu }^{\left( 2\right) }\left(
r;q\right) $. A similar relation was earlier found by Rahman \cite{Ra}. This
transformation shows that our basic sine $S\left( \eta ;\omega \right) $ and
basic cosine $C\left( \eta ;\omega \right) $ functions are just multiples of 
$J_{1/2}^{\left( 2\right) }\left( 2\omega ;q\right) $ and $J_{-1/2}^{\left(
2\right) }\left( 2\omega ;q\right) $, namely, 
\begin{equation}
S\left( \eta ;\omega \right) =\frac{\left( q;q\right) _{\infty }}{\left(
q^{1/2};q\right) _{\infty }}\ \frac{\omega ^{1/2}}{\left( -q\omega
^{2};q^{2}\right) _{\infty }}\ J_{1/2}^{\left( 2\right) }\left( 2\omega
;q\right) ,
\end{equation}
\begin{equation}
C\left( \eta ;\omega \right) =\frac{\left( q;q\right) _{\infty }}{\left(
q^{1/2};q\right) _{\infty }}\ \frac{\omega ^{1/2}}{\left( -q\omega
^{2};q^{2}\right) _{\infty }}\ J_{-1/2}^{\left( 2\right) }\left( 2\omega
;q\right) .
\end{equation}
The main properties of zeros of the $q$-Bessel functions $J_{\nu }^{\left(
2\right) }\left( r;q\right) $ were established in Ismail's papers \cite{Is:M}
and \cite{Is:Mo} by a different method. This gives independent proofs of our
Theorems 5.1--5.4. Some monotonicity properties of zeros of $J_{\nu
}^{\left( 2\right) }\left( r;q\right) $ were discussed in \cite{Is:Mu}.
Chen, Ismail, and Muttalib \cite{Ch:Is:Mu} have found a complete asymptotic
expansion of $J_{\nu }^{\left( 2\right) }\left( r;q\right) \mathrm{\ for}$
the large argument, 
\begin{eqnarray}
J_{\nu }^{\left( 2\right) }\left( r;q\right) &=&\frac{\left(
q^{1/2};q\right) _{\infty }}{2\left( q;q\right) _{\infty }}\ \left( \frac{r}{%
2}\right) ^{\nu } \\
&&\times \Biggl[\left( i\ \frac{r}{2}\ q^{\left( \nu +1/2\right)
/2};q^{1/2}\right) _{\infty }  \notag \\
&&\quad \times \sum_{n=0}^{\infty }q^{n/2}\frac{\left( q^{\nu +1/2};q\right)
_{n}}{\left( q;q\right) _{n}}\left( i\ \frac{r}{2}\ q^{\left( \nu
+1/2\right) /2};q^{1/2}\right) _{n}^{-1}  \notag \\
&&\quad +\left( -i\ \frac{r}{2}\ q^{\left( \nu +1/2\right)
/2};q^{1/2}\right) _{\infty }  \notag \\
&&\quad \times \sum_{n=0}^{\infty }q^{n/2}\frac{\left( q^{\nu +1/2};q\right)
_{n}}{\left( q;q\right) _{n}}\left( -i\ \frac{r}{2}\ q^{\left( \nu
+1/2\right) /2};q^{1/2}\right) _{n}^{-1}\Biggr].  \notag
\end{eqnarray}
This follows also from Exersises 3.15 and 3.8 of \cite{Ga:Ra}. Equations
(5.34)--(5.36) result in (5.11) and (5.15).}
\end{remark}

\section{Evaluation of Some Constants\textrm{\ }}

In this section we shall find explicitly the values of the normalization
constants in the right sides of the orthogonality relations (3.26)--(3.27)
for the basic sine and basic cosine functions. First, we evaluate the
integral 
\begin{eqnarray}
2k\left( \omega \right) &=&\int_{0}^{\pi }\left( C^{2}\left( \cos \theta
;\omega \right) +S^{2}\left( \cos \theta ;\omega \right) \right) \  \\
&&\ \ \times \frac{\ \left( e^{2i\theta },\ e^{-2i\theta };\ q\right)
_{\infty }}{\left( q^{1/2}e^{2i\theta },\ q^{1/2}e^{-2i\theta };\ q\right)
_{\infty }}\ d\theta  \notag \\
&=&\int_{0}^{\pi }C\left( \cos \theta ,-\cos \theta ;\omega \right) \ \frac{%
\ \left( e^{2i\theta },\ e^{-2i\theta };\ q\right) _{\infty }}{\left(
q^{1/2}e^{2i\theta },\ q^{1/2}e^{-2i\theta };\ q\right) _{\infty }}\ d\theta
,  \notag
\end{eqnarray}
where we have used the identity (4.14) of \cite{Su}, 
\begin{equation}
C\left( x,-x;\omega \right) =C^{2}\left( x;\omega \right) +S^{2}\left(
x;\omega \right) .
\end{equation}
In view of (2.3), for $|\omega |<1$ one can write 
\begin{eqnarray}
2\ \frac{\left( -q\omega ^{2};q^{2}\right) _{\infty }}{\left( -\omega
^{2};q^{2}\right) _{\infty }}\ k\left( \omega \right) &=&\sum_{n=0}^{\infty
}\left( -\omega ^{2}\right) ^{n}\frac{\left( q^{1/2};q\right) _{n}}{\left(
q,-q,-q^{1/2};q\right) _{n}} \\
&&\times \int_{0}^{\pi }\frac{\ \left( e^{2i\theta },\ e^{-2i\theta };\
q\right) _{\infty }}{\left( q^{n+1/2}e^{2i\theta },\ q^{n+1/2}e^{-2i\theta
};\ q\right) _{\infty }}\ d\theta .  \notag
\end{eqnarray}
The integral in the right side is a special case of the Askey--Wilson
integral \cite{As:Wi}, 
\begin{eqnarray}
&&\int_{0}^{\pi }\frac{\ \left( e^{2i\theta },\ e^{-2i\theta };\ q\right)
_{\infty }}{\left( q^{n+1/2}e^{2i\theta },\ q^{n+1/2}e^{-2i\theta };\
q\right) _{\infty }}\ d\theta \\
&&\ =\frac{2\pi \left( q^{2n+2};q\right) _{\infty }}{\left(
q,-q^{n+1/2},q^{n+1},-q^{n+1},q^{n+1},-q^{n+1},-q^{n+3/2};q\right) _{\infty }%
}.  \notag
\end{eqnarray}
Therefore, 
\begin{eqnarray}
2\ \frac{\left( -q\omega ^{2};q^{2}\right) _{\infty }}{\left( -\omega
^{2};q^{2}\right) _{\infty }}\ k\left( \omega \right) &=&\frac{2\pi \left(
q^{1/2};q\right) _{\infty }}{\left( q,q,-q,-q^{1/2};q\right) _{\infty }} \\
&&\times \sum_{n=0}^{\infty }\frac{\left( -\omega ^{2}\right) ^{n}}{%
1-q^{n+1/2}},  \notag
\end{eqnarray}
where we have used the identity 
\begin{equation*}
\left( q^{2n+2};q\right) _{\infty }=\left(
q^{n+1},-q^{n+1},q^{n+3/2},-q^{n+3/2};q\right) _{\infty }.
\end{equation*}
But, 
\begin{eqnarray*}
&&\sum_{n=0}^{\infty }\frac{\left( -\omega ^{2}\right) ^{n}}{1-q^{n+1/2}} \\
\ &&\ =\frac{1}{1-q^{1/2}}\ \ _{2}\varphi _{1}\left( 
\begin{array}{c}
q,\ q^{1/2} \\ 
q^{3/2}
\end{array}
;\ q,-\omega ^{2}\right) \\
&&\ =\frac{\left( q,-q^{1/2}\omega ^{2};q\right) _{\infty }}{\left(
q^{1/2},-\omega ^{2};q\right) _{\infty }}\ \ _{2}\varphi _{1}\left( 
\begin{array}{c}
\ q^{1/2},\ -\omega ^{2} \\ 
-q^{1/2}\omega ^{2}
\end{array}
;\ q,q\right)
\end{eqnarray*}
by (III.1) of \cite{Ga:Ra}. The last line provides an analytic continuation
of this sum in the complex $\omega $-plane. Finally, we obtain 
\begin{eqnarray}
\ k\left( \omega \right) &=&\frac{1}{2}\int_{0}^{\pi }\left( C^{2}\left(
\cos \theta ;\omega \right) +S^{2}\left( \cos \theta ;\omega \right) \right)
\\
&&\ \ \times \frac{\ \left( e^{2i\theta },\ e^{-2i\theta };\ q\right)
_{\infty }}{\left( q^{1/2}e^{2i\theta },\ q^{1/2}e^{-2i\theta };\ q\right)
_{\infty }}\ d\theta  \notag \\
&=&\pi \frac{\left( q^{1/2},-q^{1/2}\omega ^{2};q\right) _{\infty }}{\left(
q,-\omega ^{2};q\right) _{\infty }}\ \frac{\left( -\omega ^{2};q^{2}\right)
_{\infty }}{\left( -q\omega ^{2};q^{2}\right) _{\infty }} \\
&&\ \ \times \ _{2}\varphi _{1}\left( 
\begin{array}{c}
\ q^{1/2},\ -\omega ^{2} \\ 
-q^{1/2}\omega ^{2}
\end{array}
;\ q,q\right) .  \notag
\end{eqnarray}
The second line gives the large asymptotic of the function $k\left( \omega
\right) $, 
\begin{equation}
k\left( \omega \right) =\pi \frac{\left( -q^{1/2}\omega ^{2};q\right)
_{\infty }}{\left( -\omega ^{2};q\right) _{\infty }}\ \frac{\left( -\omega
^{2};q^{2}\right) _{\infty }}{\left( -q\omega ^{2};q^{2}\right) _{\infty }}\
\left[ 1+\text{o}\left( \omega ^{-2}\right) \right] ,
\end{equation}
as $\omega ^{2}\rightarrow \infty $ but $\omega ^{2}\neq -q^{-n-1/2}$ for a
positive integer $n$. In particular, when $\omega =\omega _{n}$ are large
zeros of the basic sine function $S\left( \eta ;\omega \right) $ one gets as 
$n\rightarrow \infty $%
\begin{equation}
k\left( \omega _{n}\right) \sim k\left( q^{1/4-n}\right) \sim 2\pi \frac{%
\left( -q;q\right) _{\infty }^{2}}{\left( -q^{1/2};q\right) _{\infty }^{2}}
\end{equation}
by (5.11) and (I.9) of \cite{Ga:Ra}.

With the aid of (6.6)--(6.7) one can now rewrite (3.31) and (3.32) in more
explicit form, 
\begin{eqnarray}
&&\int_{0}^{\pi }C^{2}\left( \cos \theta ;\omega \right) \ \frac{\ \left(
e^{2i\theta },\ e^{-2i\theta };\ q\right) _{\infty }}{\left(
q^{1/2}e^{2i\theta },\ q^{1/2}e^{-2i\theta };\ q\right) _{\infty }}\ d\theta
\\
&&\ =\ k\left( \omega \right) +\dfrac{\pi \left( q^{1/2};\ q\right) _{\infty
}^{2}}{\omega \left( q;\ q\right) _{\infty }^{2}}\ \ C\left( \eta ;\omega
\right) S\left( \eta ;\omega \right)  \notag
\end{eqnarray}
and 
\begin{eqnarray}
&&\int_{0}^{\pi }S^{2}\left( \cos \theta ;\omega \right) \ \frac{\ \left(
e^{2i\theta },\ e^{-2i\theta };\ q\right) _{\infty }}{\left(
q^{1/2}e^{2i\theta },\ q^{1/2}e^{-2i\theta };\ q\right) _{\infty }}\ d\theta
\\
&&\ =\ k\left( \omega \right) -\dfrac{\pi \left( q^{1/2};\ q\right) _{\infty
}^{2}}{\omega \left( q;\ q\right) _{\infty }^{2}}\ \ C\left( \eta ;\omega
\right) S\left( \eta ;\omega \right) .  \notag
\end{eqnarray}
These basic integrals are, obviously, $q$-extensions of the following
elementary integrals 
\begin{equation}
\int_{-1}^{1}\cos ^{2}\omega x\ dx=1+\frac{1}{\omega }\ \sin \omega \cos
\omega ,
\end{equation}
\begin{equation}
\int_{-1}^{1}\sin ^{2}\omega x\ dx=1-\frac{1}{\omega }\ \sin \omega \cos
\omega ,
\end{equation}
respectively.

When $\omega $ satisfies the boundary condition (3.32) the last terms in the
right sides of (6.10) and (6.11) vanish and we obtain the values of the
normalization constants in the orthogonality relations (3.26)--(3.28) in
terms of the function $k\left( \omega \right) $ defined by (6.7).

\section{Orthogonality Relations for $q$-Exponential Functions}

Euler's formula, 
\begin{equation}
e^{i\omega x}=\cos \omega x+i\sin \omega x,
\end{equation}
allows us to rewrite the orthogonality relations for the trigonometric
functions (1.6)--(1.8) in a complex form, 
\begin{equation}
\frac{1}{2l}\int_{-l}^{l}\exp \left( i\frac{\pi m}{l}x\right) \ \exp \left(
-i\frac{\pi n}{l}x\right) \ dx=\delta _{mn},
\end{equation}
where 
\begin{equation}
\delta _{mn}=\left\{ 
\begin{array}{ccc}
1 & \text{if} & m=n, \\ 
0 & \text{if} & m\neq n.
\end{array}
\right.
\end{equation}

The $q$-analog of Euler's formula (7.1) is 
\begin{equation}
\mathcal{E}_{q}\left( x;i\omega \right) =C_{q}\left( x;\omega \right)
+iS_{q}\left( x;\omega \right) ,
\end{equation}
where $\mathcal{E}_{q}\left( x;\alpha \right) $ with $\alpha =i\omega $ is
the $q$-exponential function introduced in \cite{Is:Zh} (see also \cite
{At:Su} and \cite{Su}, we shall use the same notations as in \cite{Su}); $%
C_{q}\left( x;\omega \right) $ and $S_{q}\left( x;\omega \right) $ are basic
cosine and sine functions defined by (2.1) and (2.2), respectively. Our
orthogonality relations for the basic trigonometric functions (3.26)--(3.28)
result in the following \textit{orthogonality property }for the $q$%
-exponential function 
\begin{eqnarray}
&&\frac{1}{2k\left( \omega _{n}\right) }\int_{0}^{\pi }\mathcal{E}_{q}\left(
\cos \theta ;i\omega _{m}\right) \ \mathcal{E}_{q}\left( \cos \theta
;-i\omega _{n}\right) \\
&&\qquad \quad \times \ \frac{\ \left( e^{2i\theta },\ e^{-2i\theta };\
q\right) _{\infty }}{\left( q^{1/2}e^{2i\theta },\ q^{1/2}e^{-2i\theta };\
q\right) _{\infty }}\ d\theta =\delta _{mn},  \notag
\end{eqnarray}
where $\omega _{m},\omega _{n}=0,\pm \omega _{1},\pm \omega _{2},\pm \omega
_{3},...$ and $\omega _{0}=0,\omega _{1},\omega _{2},\omega _{3},...$, are
nonnegative zeros of the basic sine function $S\left( \eta ;\omega \right) $
arranged in ascending order of magnitude; the normalization constants $%
k\left( \omega _{n}\right) $ are defined by (6.7).

A basic analog of 
\begin{equation}
e^{i\omega \left( x+y\right) }=\cos \omega \left( x+y\right) +i\sin \omega
\left( x+y\right)
\end{equation}
is 
\begin{equation}
\mathcal{E}_{q}\left( x,y;i\omega \right) =C_{q}\left( x,y;\omega \right)
+iS_{q}\left( x,y;\omega \right) ,
\end{equation}
see \cite{Is:Zh} and \cite{Su}. The general exponential function on a $q$%
-quadratic grid $\mathcal{E}_{q}\left( x,y;i\omega \right) $ has the
following orthogonality property.

\begin{theorem}
\begin{eqnarray}
&&\int_{0}^{\pi }\mathcal{E}_{q}\left( \cos \theta ,\cos \varphi ;i\omega
_{m}\right) \ \mathcal{E}_{q}\left( \cos \theta ,\cos \varphi ^{\prime
};-i\omega _{n}\right) \\
&&\qquad \times \ \frac{\ \left( e^{2i\theta },\ e^{-2i\theta };\ q\right)
_{\infty }}{\left( q^{1/2}e^{2i\theta },\ q^{1/2}e^{-2i\theta };\ q\right)
_{\infty }}\ d\theta  \notag \\
&&\qquad =2k\left( \omega _{n}\right) \ \mathcal{E}_{q}\left( \cos \varphi
;i\omega _{n}\right) \ \mathcal{E}_{q}\left( \cos \varphi ^{\prime
};-i\omega _{n}\right) \ \delta _{mn},  \notag
\end{eqnarray}
where $\omega _{m},\omega _{n}=0,\pm \omega _{1},\pm \omega _{2},\pm \omega
_{3},...$ and $\omega _{0}=0,\omega _{1},\omega _{2},\omega _{3},...$, are
nonnegative zeros of the basic sine function $S\left( \eta ;\omega \right) $
arranged in ascending order of magnitude; the normalization constants $%
k\left( \omega _{n}\right) $ are defined by \textrm{(6.7)}.
\end{theorem}

\proof%
Using of the ``addition'' theorem for basic exponential functions \cite{Su}, 
\begin{equation}
\mathcal{E}_{q}\left( x,y;i\omega \right) =\mathcal{E}_{q}\left( x;i\omega
\right) \ \mathcal{E}_{q}\left( y;i\omega \right) ,
\end{equation}
and the orthogonality relation (7.5) one gets (7.8). 
\endproof%

In a similar fashion, we can establish the following results.

\begin{theorem}
\begin{eqnarray}
&&\int_{0}^{\pi }C\left( \cos \theta ,\cos \varphi ;\omega _{m}\right) \
C\left( \cos \theta ,\cos \varphi ^{\prime };\omega _{n}\right) \\
&&\qquad \times \ \frac{\ \left( e^{2i\theta },\ e^{-2i\theta };\ q\right)
_{\infty }}{\left( q^{1/2}e^{2i\theta },\ q^{1/2}e^{-2i\theta };\ q\right)
_{\infty }}\ d\theta  \notag \\
&\ &\ =\left\{ 
\begin{tabular}{lll}
$0$ & if & $m\neq n,$ \\ 
$k\left( \omega _{n}\right) \ C\left( \cos \varphi ,-\cos \varphi ^{\prime
};\omega _{n}\right) $ & if & $m=n;$%
\end{tabular}
\right.  \notag
\end{eqnarray}
\begin{eqnarray}
&&\int_{0}^{\pi }S\left( \cos \theta ,\cos \varphi ;\omega _{m}\right) \
S\left( \cos \theta ,\cos \varphi ^{\prime };\omega _{n}\right) \\
&&\qquad \times \ \frac{\ \left( e^{2i\theta },\ e^{-2i\theta };\ q\right)
_{\infty }}{\left( q^{1/2}e^{2i\theta },\ q^{1/2}e^{-2i\theta };\ q\right)
_{\infty }}\ d\theta  \notag \\
&\ &\ =\left\{ 
\begin{tabular}{lll}
$0$ & if & $m\neq n,$ \\ 
$k\left( \omega _{n}\right) \ C\left( \cos \varphi ,-\cos \varphi ^{\prime
};\omega _{n}\right) $ & if & $m=n;$%
\end{tabular}
\right.  \notag
\end{eqnarray}
and 
\begin{eqnarray}
&&\int_{0}^{\pi }C\left( \cos \theta ,\cos \varphi ;\omega _{m}\right) \
S\left( \cos \theta ,\cos \varphi ^{\prime };\omega _{n}\right) \\
&&\qquad \times \ \frac{\ \left( e^{2i\theta },\ e^{-2i\theta };\ q\right)
_{\infty }}{\left( q^{1/2}e^{2i\theta },\ q^{1/2}e^{-2i\theta };\ q\right)
_{\infty }}\ d\theta  \notag \\
\ \ &&\ =\left\{ 
\begin{tabular}{lll}
$0$ & if & $m\neq n,$ \\ 
$k\left( \omega _{n}\right) \ S\left( \cos \varphi ,-\cos \varphi ^{\prime
};\omega _{n}\right) $ & if & $m=n;$%
\end{tabular}
\right.  \notag
\end{eqnarray}
where $\omega _{m},\omega _{n}=\omega _{1},\omega _{2},\omega _{3},...$, are
positive zeros of the basic sine function $S\left( \eta ;\omega \right) $
arranged in ascending order of magnitude; the normalization constants $%
k\left( \omega _{n}\right) $ are defined by \textrm{(6.7)}.
\end{theorem}

\proof%
Use the ``addition'' theorem for the basic trigonometric functions \cite{Su}
and the orthogonality relations (3.26)--(3.28). 
\endproof%

\section{Basic Fourier Series}

By analogy with (1.2) we can now introduce a $q$-version of Fourier series, 
\begin{equation}
f\left( \cos \theta \right) =a_{0}+\sum_{n=1}^{\infty }\left(
a_{n}C_{q}\left( \cos \theta ;\omega _{n}\right) +b_{n}S_{q}\left( \cos
\theta ;\omega _{n}\right) \right) ,
\end{equation}
where $\omega _{0}=0,\omega _{1},\omega _{2},\omega _{3},...$, are
nonnegative zeros of the basic sine function $S\left( \eta ;\omega \right) $
arranged in ascending order of magnitude, and 
\begin{eqnarray}
a_{0} &=&\frac{1}{2k\left( 0\right) }\ \int_{0}^{\pi }f\left( \cos \theta
\right) \\
&&\qquad \times \ \frac{\ \left( e^{2i\theta },\ e^{-2i\theta };\ q\right)
_{\infty }}{\left( q^{1/2}e^{2i\theta },\ q^{1/2}e^{-2i\theta };\ q\right)
_{\infty }}\ d\theta ,  \notag
\end{eqnarray}
\begin{eqnarray}
a_{n} &=&\frac{1}{k\left( \omega _{n}\right) }\ \int_{0}^{\pi }f\left( \cos
\theta \right) \ C_{q}\left( \cos \theta ;\omega _{n}\right) \\
&&\qquad \times \ \frac{\ \left( e^{2i\theta },\ e^{-2i\theta };\ q\right)
_{\infty }}{\left( q^{1/2}e^{2i\theta },\ q^{1/2}e^{-2i\theta };\ q\right)
_{\infty }}\ d\theta ,  \notag
\end{eqnarray}
\begin{eqnarray}
b_{n} &=&\frac{1}{k\left( \omega _{n}\right) }\ \int_{0}^{\pi }f\left( \cos
\theta \right) \ S_{q}\left( \cos \theta ;\omega _{n}\right) \\
&&\qquad \times \ \frac{\ \left( e^{2i\theta },\ e^{-2i\theta };\ q\right)
_{\infty }}{\left( q^{1/2}e^{2i\theta },\ q^{1/2}e^{-2i\theta };\ q\right)
_{\infty }}\ d\theta .  \notag
\end{eqnarray}
The complex form of the basic Fourier series (8.1) is 
\begin{equation}
f\left( \cos \theta \right) =\sum_{n=-\infty }^{\infty }c_{n}\ \mathcal{E}%
_{q}\left( \cos \theta ;i\omega _{n}\right)
\end{equation}
with 
\begin{eqnarray}
c_{n} &=&\frac{1}{2k\left( \omega _{n}\right) }\ \int_{0}^{\pi }f\left( \cos
\theta \right) \ \mathcal{E}_{q}\left( \cos \theta ;-i\omega _{n}\right) \\
&&\qquad \quad \times \ \frac{\ \left( e^{2i\theta },\ e^{-2i\theta };\
q\right) _{\infty }}{\left( q^{1/2}e^{2i\theta },\ q^{1/2}e^{-2i\theta };\
q\right) _{\infty }}\ d\theta  \notag
\end{eqnarray}
where $\omega _{n}=0,\pm \omega _{1},\pm \omega _{2},\pm \omega _{3},...$
and $\omega _{0}=0,\omega _{1},\omega _{2},\omega _{3},...$, are nonnegative
zeros of the basic sine function $S\left( \eta ;\omega \right) $ arranged in
ascending order of magnitude; the normalization constants $k\left( \omega
_{n}\right) $ are defined by (6.7). These expressions, of course, merely
indicate how the coefficients of our basic Fourier series are to be
determined on the hypothesis that the expansion exists and is uniformly
convergent. We shall study the question of convergence of the series (8.1)
and (8.5) in the next sections.

The $q$-Fourier series of $f$ in either of the forms (8.1) and (8.5) will be
denoted in a usual manner by $\mathbf{S}\left[ f\right] $.

\section{Completeness of the $q$-Trigonometric System}

Completeness of the trigonometric system $\left\{ e^{i\pi nx}\right\}
_{n=-\infty }^{\infty }$ on the interval $\left( -1,1\right) $ is one of the
fundamental facts in the theory of trigonometric series (see, for example, 
\cite{AKh}, \cite{Ko:Fo}, \cite{Le}, \cite{Lev} and \cite{Zyg}). In this
section we shall prove a similar property for the system of basic
trigonometric function $\left\{ \mathcal{E}_{q}\left( x;i\omega _{n}\right)
\right\} $, where $\omega _{n}=0,\pm \omega _{1},\pm \omega _{2},\pm \omega
_{3},...$ and $\omega _{0}=0,\omega _{1},\omega _{2},\omega _{3},...$, are
nonnegative zeros of the basic sine function $S\left( \eta ;\omega \right) $
arranged in ascending order of magnitude. But first we need to discuss
connections between the basic trigonometric functions and the continuous $q$%
-Hermite polynomials.

The continuous $q$-Hermite polynomials, 
\begin{equation}
H_{n}\left( \cos \theta |q\right) =\sum_{k=0}^{n}\frac{\left( q;q\right) _{n}%
}{\left( q;q\right) _{k}\left( q;q\right) _{n-k}}\ e^{i\left( n-2k\right)
\theta },
\end{equation}
have two generating functions, 
\begin{equation}
\sum_{n=0}^{\infty }\frac{r^{n}}{\left( q;q\right) _{n}}\ H_{n}\left( \cos
\theta |q\right) =\frac{1}{\left( re^{i\theta },re^{-i\theta };q\right)
_{\infty }},
\end{equation}
when $|r|<1$ and 
\begin{equation}
\sum_{n=0}^{\infty }\ \frac{q^{n^{2}/4}}{\left( q;q\right) _{n}}\ \alpha
^{n}H_{n}\left( \cos \theta |q\right) =\left( q\alpha ^{2};q^{2}\right)
_{\infty }\mathcal{E}_{q}\left( \cos \theta ;\alpha \right)
\end{equation}
(see, for example, \cite{Ga:Ra}, \cite{Is:Zh}, and \cite{Su}).

\begin{lemma}
The following functions 
\begin{eqnarray}
e\left( x,\alpha \right) &=&\left( q\alpha ^{2};q^{2}\right) _{\infty }%
\mathcal{E}_{q}\left( x;\alpha \right) , \\
s(\omega ) &=&\left( -q\omega ^{2};q^{2}\right) \ S\left( \eta ;\omega
\right) \\
&=&\frac{1}{2i}\ \left( e(x,i\omega )-e(x,-i\omega )\right) ,  \notag
\end{eqnarray}
and 
\begin{eqnarray}
c(\omega ) &=&\left( -q\omega ^{2};q^{2}\right) \ C\left( \eta ;\omega
\right) \\
&=&\frac{1}{2}\ \left( e(x,i\omega )+e(x,-i\omega )\right)  \notag
\end{eqnarray}
are entire functions in $\alpha $ and $\omega $, respectively, of order zero
for all real values of $x$.
\end{lemma}

\proof%
The generating function (9.3) gives a power series expansion for the
function (9.4),

\begin{equation}
e\left( x,\alpha \right) =\sum_{n=0}^{\infty }h_{n}\ \alpha ^{n}
\end{equation}
with 
\begin{equation}
h_{n}=h_{n}(x)=\frac{q^{n^{2}/4}}{\left( q;q\right) _{n}}\ H_{n}\left(
x|q\right) .
\end{equation}
The radius of convergence of this series is infinity, because 
\begin{eqnarray}
\frac{1}{R} &=&\lim_{n\rightarrow \infty }\left( \left| h_{n}\right| \right)
^{1/n} \\
&=&\lim_{n\rightarrow \infty }\left( \left| \frac{q^{n^{2}/4}}{\left(
q;q\right) _{n}}\ H_{n}\left( x|q\right) \right| \right) ^{1/n}=0.  \notag
\end{eqnarray}
Thus, $e\left( x,\alpha \right) $ is an entire function in $\alpha $. The
order of this entire function is \cite{Lev} 
\begin{eqnarray}
&&\lim_{n\rightarrow \infty }\left( \frac{n\log n}{-\log \left| h_{n}\right| 
}\right) \\
&&\ =\lim_{n\rightarrow \infty }\left( \frac{n\log n}{-\log \left|
q^{n^{2}/4}\ H_{n}\left( x|q\right) /\left( q;q\right) _{n}\right| }\right)
=0.  \notag
\end{eqnarray}
Functions (9.5) and (9.6) are just a sum or difference of two functions of
type (9.4), so they are also entire functions of order zero. This proves the
lemma. 
\endproof%

The next step is to establish the following inequalities.

\begin{lemma}
Let $-\cosh \tau \leq -1\leq x\leq 1\leq \cosh \tau ,$ where $x=\cos \theta $%
, $0\leq \theta \leq \pi $, and $\tau \geq 0$. Then 
\begin{equation}
\left| e\left( \cos \theta ;\alpha \right) \right| \leq e\left( \cosh \tau
;\left| \alpha \right| \right)
\end{equation}
and 
\begin{equation}
e\left( \cosh \tau ;\left| \alpha \right| \right) \leq e\left( \cosh \tau
_{1};\left| \alpha \right| \right)
\end{equation}
if $\tau <\tau _{1}$.
\end{lemma}

\proof%
One can rewrite (9.1) as 
\begin{equation}
H_{n}\left( \cos \theta |q\right) =2\sum_{k=0}^{\left[ n/2\right] }\frac{%
\left( q;q\right) _{n}}{\left( q;q\right) _{k}\left( q;q\right) _{n-k}}\
\cos \left( n-2k\right) \theta .
\end{equation}
Thus, 
\begin{eqnarray}
\left| H_{n}\left( \cos \theta |q\right) \right| &\leq &2\sum_{k=0}^{\left[
n/2\right] }\frac{\left( q;q\right) _{n}}{\left( q;q\right) _{k}\left(
q;q\right) _{n-k}}\ \cosh \left( n-2k\right) \tau \\
&=&H_{n}\left( \cosh \tau |q\right) .  \notag
\end{eqnarray}
Estimating both sides of (9.3) gives 
\begin{eqnarray*}
&&\left| \left( q\alpha ^{2};q^{2}\right) _{\infty }\mathcal{E}_{q}\left(
\cos \theta ;\alpha \right) \right| \\
&&\ \leq \sum_{n=0}^{\infty }\ \frac{q^{n^{2}/4}}{\left( q;q\right) _{n}}\
\left| \alpha \right| ^{n}\left| H_{n}\left( \cos \theta |q\right) \right| \\
&&\ \leq \sum_{n=0}^{\infty }\ \frac{q^{n^{2}/4}}{\left( q;q\right) _{n}}\
\left| \alpha \right| ^{n}H_{n}\left( \cosh \tau |q\right) \\
&&\ =\left( q\left| \alpha \right| ^{2};q^{2}\right) _{\infty }\mathcal{E}%
_{q}\left( \cosh \tau ;\left| \alpha \right| \right)
\end{eqnarray*}
by (9.14) and (9.3). This proves (9.11). The monotonicity property (9.12)
follows from the monotonicity of the hyperbolic cosine function. 
\endproof%

It is clear
that the system $\left\{ \mathcal{E}_{q}\left( x;i\omega _{n}\right)
\right\} _{n=-\infty }^{\infty }$ is complete if the equivalent system $%
\left\{ e\left( x,i\omega _{n}\right) \right\} _{n=-\infty }^{\infty }$ is
closed.

Suppose that the system $\left\{ e\left( x,i\omega _{n}\right) \right\}
_{n=-\infty }^{\infty }$ is not closed on $\left( -1,1\right) $. This means
that there exists at least one function $\phi (x)$, not identically zero,
such that 
\begin{equation}
\int_{-1}^{1}\phi (x)\ e\left( x,i\omega _{n}\right) \ \rho (x)\ dx=0,\qquad
n=0,\pm 1,\pm 2,...,
\end{equation}
where $\rho (x)$ is the absolutely continuous measure in the orthogonality
relation (7.5). Then, the function 
\begin{equation}
f(\omega )=\int_{-1}^{1}\phi (x)\ e\left( x,i\omega _{n}\right) \ \rho (x)\
dx
\end{equation}
is an entire function of order zero and $f(\omega _{n})=0$ for all $n=0,\pm
1,\pm 2,...$. Thus the study of closure amounts to the study of zeros of a
certain entire function.

Suppose that $\phi (x)\ $is integrable on $\left( -1,1\right) $, 
\begin{equation}
\int_{-1}^{1}\left| \phi (x)\right| \ \rho (x)\ dx=A<\infty .
\end{equation}
Then 
\begin{eqnarray}
\left| f(\omega )\right| &\leq &\int_{-1}^{1}\left| \phi (x)\ e\left(
x,i\omega \right) \right| \ \rho (x)\ dx \\
&\leq &e\left( \cosh \tau ,\left| \omega \right| \right) \   \notag \\
&&\times \int_{-1}^{1}\left| \phi (x)\right| \ \rho (x)\ dx  \notag \\
&=&A\ e\left( \cosh \tau ,\left| \omega \right| \right)  \notag
\end{eqnarray}
by (9.11) and (9.17).

Consider the quotient 
\begin{equation}
g(\omega )=\frac{f(\omega )}{s(\omega )}
\end{equation}
of two entire functions, $f(\omega )$ and $s(\omega )$ defined by (9.16) and
(9.5), respestively. The functions $f(\omega )$ and $s(\omega )$ have the
same zeros, so $g(\omega )$ is an entire function. The order of this entire
function is zero because both $f(\omega )$ and $s(\omega )$ are of order
zero (see \cite{Lev}, Corollary of Theorem 12 on p. 24). Moreover, this
function $g(\omega )$ is bounded on a straight line parallel to the
imaginary axis. Indeed, let $\omega =\gamma +i\delta $. Using the same
arguments as in Section 5 one can see that 
\begin{equation}
\lim_{\left| \delta \right| \rightarrow \infty }\left| \frac{s(i\delta )}{%
e(\eta ,|\delta |)}\right| <\infty .
\end{equation}
From this condition and the inequality (9.18), it follows that the entire
function $g(\omega )$ is bounded on the imaginary axis. But an entire
function of order zero bounded on a line must be a constant (see Theorems
21--22 and Corollary on pp. 49--51 of \cite{Lev}). Then, 
\begin{equation}
f(\omega )=c\ s(\omega )
\end{equation}
and, therefore, 
\begin{eqnarray}
|c| &=&\left| \int_{-1}^{1}\phi (x)\ \frac{e\left( x,i\omega \right) }{%
s(\omega )}\ \rho (x)\ dx\right| \\
&\leq &\int_{-1}^{1}\left| \phi (x)\ \frac{e\left( x,i\omega \right) }{%
s(\omega )}\right| \ \rho (x)\ dx  \notag \\
&\leq &A\ \left| \frac{\mathcal{E}_{q}\left( \cosh \tau ;|\omega |\right) }{%
S(\cosh \tau _{1};|\omega |)}\right| \rightarrow 0  \notag
\end{eqnarray}
as $\left| \omega \right| \rightarrow \infty $ and $\tau <\tau _{1}$. Thus, $%
f(\omega )$ is identically zero and the function $\phi (x)$ does not exist.

We have established the following theorem.

\begin{theorem}
The system of the basic trigonometric function $\left\{ \mathcal{E}%
_{q}\left( x;i\omega _{n}\right) \right\} $, where $n=0,\pm 1,\pm 2,...$ and 
$\omega _{0}=0,\omega _{1},\omega _{2},\omega _{3},...$, are nonnegative
zeros of the basic sine function $S\left( \eta ;\omega \right) $ arranged in
ascending order of magnitude, is complete on $\left( -1,1\right) $.
\end{theorem}

As corollaries we have the following results.

\begin{theorem}
If $f(x)$ and $g(x)$ have the same $q$-Fourier series, then $f\equiv g$.
\end{theorem}

\proof%
The $q$-Fourier coefficients of $f-g$ all vanish, so that $f-g\equiv 0$. 
\endproof%

\begin{theorem}
If $f(x)$ is continuous and $\mathbf{S}\left[ f\right] $, the $q$-Fourier
series of function $f$, converges uniformly, then its sum is $f(x)$.
\end{theorem}

\proof%
Let $g(x)$ denote the sum of $\mathbf{S}\left[ f\right] $, the $q$-Fourier
series in the right side of (8.5). Then the coefficients of $\mathbf{S}%
\left[ f\right] $ are $q$-Fourier coefficients of $g$. Hence, $\mathbf{S}%
\left[ f\right] =\mathbf{S}\left[ g\right] $, so that $f\equiv g$ and, $f$
and $g$ being continuous, $f(x)\equiv g(x)$. 
\endproof%

Bessel's inequality for the $q$-trigonometric system $\left\{ \mathcal{E}%
_{q}\left( x;i\omega _{n}\right) \right\} _{n=-\infty }^{\infty }$, where $%
\omega _{0}=0,\omega _{1},\omega _{2},\omega _{3},...$, are nonnegative
zeros of the basic sine function $S\left( \eta ;\omega \right) $ arranged in
ascending order of magnitude, takes the form 
\begin{equation}
\sum_{n=-N}^{N}\left| c_{n}\right| ^{2}\leq \int_{-1}^{1}\left| f(x)\right|
^{2}\ \rho (x)\ dx
\end{equation}
provided $f\in L_{\rho }^{2}\left( -1,1\right) $, which means that $\left|
f(x)\right| ^{2}$ is integrable on $\left( -1,1\right) $ with respect to the
weight function $\rho (x)$ in the orthogonality relation (7.5). Here $c_{n}$
are the $q$-Fourier coefficients of $f(x)$ defined by (8.6). When $%
N\rightarrow \infty $ we get Parseval's formula 
\begin{equation}
\sum_{n=-\infty }^{\infty }\left| c_{n}\right| ^{2}=\int_{-1}^{1}\left|
f(x)\right| ^{2}\ \rho (x)\ dx
\end{equation}
due to the completeness of the $q$-trigonometric system $\left\{ \mathcal{E}%
_{q}\left( x;i\omega _{n}\right) \right\} _{n=-\infty }^{\infty }$ and the
space $L_{\rho }^{2}\left( -1,1\right) $ \cite{AKh}, \cite{Ko:Fo}. It
follows that the $q$-Fourier coefficients $c_{n}$ tend to zero if $f\in
L_{\rho }^{2}\left( -1,1\right) $.

\section{Bilinear Generating Function}

In this section we shall derive the following bilinear generating relation, 
\begin{eqnarray}
&&\sum_{n=-\infty }^{\infty }\frac{\left( -qr^{2}\omega
_{n}^{2};q^{2}\right) _{\infty }}{\left( -q\omega _{n}^{2};q^{2}\right)
_{\infty }}\ k^{-1}\left( \omega _{n}\right) \\
&&\ \times \ \mathcal{E}_{q}\left( \cos \theta ;i\omega _{n}\right) \mathcal{%
E}_{q}\left( \cos \varphi ;-ir\omega _{n}\right)  \notag \\
&&\ =\frac{\left( q,r^{2},q^{1/2}e^{2i\theta },\ q^{1/2}e^{-2i\theta };\
q\right) _{\infty }}{\pi \left( re^{i\theta +i\varphi },re^{i\theta
-i\varphi },re^{i\varphi -i\theta },re^{-i\theta -i\varphi };q\right)
_{\infty }},  \notag
\end{eqnarray}
for the basic exponential functions. Here as before $\omega _{n}=0,\pm
\omega _{1},\pm \omega _{2},\pm \omega _{3},...$ and $\omega _{0}=0,\omega
_{1},\omega _{2},\omega _{3},...$ are nonnegative zeros of the basic sine
function $S\left( \eta ;\omega \right) $ arranged in ascending order of
magnitude. We shall use this generating function for a further investigation
of the convergence of the basic Fourier series (8.5) in the subsequent
section.

Let us establish a connecting relation of the form, 
\begin{eqnarray}
&&\frac{\left( q\alpha ^{2}r^{2};q^{2}\right) _{\infty }}{\left( q\alpha
^{2};q^{2}\right) _{\infty }}\ \mathcal{E}_{q}\left( \cos \theta ;\alpha
r\right) \\
&&\ =\frac{1}{2\pi }\int_{0}^{\pi }\frac{\left( q,r^{2},e^{2i\varphi },\
e^{-2i\varphi };\ q\right) _{\infty }}{\left( re^{i\theta +i\varphi
},re^{i\theta -i\varphi },re^{i\varphi -i\theta },re^{-i\theta -i\varphi
};q\right) _{\infty }}  \notag \\
&&\qquad \qquad \times \ \mathcal{E}_{q}\left( \cos \varphi ;\alpha \right)
\ d\varphi ,  \notag
\end{eqnarray}
where $|r|<1$. One can easily see that if we could prove the uniform
convergence in the variable $x=\cos \theta $ of the series in the left side
of (10.1), than the integral in (10.2) gives the correct values of the basic
Fourier coefficients (see (8.5)--(8.6)), which verifies the generating
relation (10.1) by Theorem 9.5. So, one needs to give a prove of (10.2)
first.

The continuous $q$-Hermite polynomials have the following bilinear
generating function (the Poisson kernel), 
\begin{eqnarray}
&&\sum_{n=0}^{\infty }\frac{r^{n}}{\left( q;q\right) _{n}}\ H_{n}\left( \cos
\theta |q\right) H_{n}\left( \cos \varphi |q\right) \\
&&\ =\frac{\left( r^{2};\ q\right) _{\infty }}{\left( re^{i\theta +i\varphi
},re^{i\theta -i\varphi },re^{i\varphi -i\theta },re^{-i\theta -i\varphi
};q\right) _{\infty }},  \notag
\end{eqnarray}
where $|r|<1$. The orthogonality relation for these polynomials is 
\begin{eqnarray}
&&\int_{0}^{\pi }\ H_{m}\left( \cos \theta |q\right) H_{n}\left( \cos \theta
|q\right) \left( e^{2i\theta },\ e^{-2i\theta };\ q\right) _{\infty }\
d\theta \\
&&\ =2\pi \frac{\left( q;q\right) _{n}}{\left( q;q\right) _{\infty }}\
\delta _{mn}  \notag
\end{eqnarray}
(see, for example, \cite{Ga:Ra}). Expanding $\mathcal{E}_{q}\left( \cos
\varphi ;\alpha \right) $ in the right side of (10.2) in the uniformly
convergent series of the $q$-Hermite polynomials with the aid of (9.3), we
get 
\begin{eqnarray}
&&\frac{1}{2\pi }\int_{0}^{\pi }\frac{\left( q,r^{2},e^{2i\varphi },\
e^{-2i\varphi };\ q\right) _{\infty }}{\left( re^{i\theta +i\varphi
},re^{i\theta -i\varphi },re^{i\varphi -i\theta },re^{-i\theta -i\varphi
};q\right) _{\infty }} \\
&&\qquad \times \left( q\alpha ^{2};q^{2}\right) _{\infty }\ \mathcal{E}%
_{q}\left( \cos \varphi ;\alpha \right) \ d\varphi  \notag \\
&&\quad =\sum_{n=0}^{\infty }\frac{q^{n^{2}/4}}{\left( q;q\right) _{n}}\
\alpha ^{n}  \notag \\
&&\qquad \times \ \frac{1}{2\pi }\int_{0}^{\pi }\frac{\left(
q,r^{2},e^{2i\varphi },\ e^{-2i\varphi };\ q\right) _{\infty }}{\left(
re^{i\theta +i\varphi },re^{i\theta -i\varphi },re^{i\varphi -i\theta
},re^{-i\theta -i\varphi };q\right) _{\infty }}  \notag \\
\qquad &&\qquad \qquad \qquad \times H_{n}\left( \cos \varphi |q\right) \
d\varphi .  \notag
\end{eqnarray}
The series in (10.3) converges uniformly when $\left| r\right| <1$. Then,
using (10.4), 
\begin{eqnarray}
&&\frac{1}{2\pi }\int_{0}^{\pi }\frac{H_{n}\left( \cos \varphi |q\right) \
\left( q,r^{2},e^{2i\varphi },\ e^{-2i\varphi };\ q\right) _{\infty }}{%
\left( re^{i\theta +i\varphi },re^{i\theta -i\varphi },re^{i\varphi -i\theta
},re^{-i\theta -i\varphi };q\right) _{\infty }}\ d\varphi \\
&&\quad =r^{n}\ H_{n}\left( \cos \theta |q\right) .  \notag
\end{eqnarray}
From (10.5), (10.6), and (9.3) we finally arrive at the connecting relation
(10.2).

Uniform convergence of the series in (10.1) can be justified with the help
of the inequality (9.11) and the corresponding asymptotic expressions. This
proves (10.1) by Theorem 9.5.

It is worth mentioning a few special cases of (10.1). When $r=0$ we obtain
the following generating function, 
\begin{eqnarray}
&&\sum_{n=-\infty }^{\infty }\frac{1}{\left( -q\omega _{n}^{2};q^{2}\right)
_{\infty }\ k\left( \omega _{n}\right) } \\
&&\ \qquad \times \ \mathcal{E}_{q}\left( \cos \theta ;i\omega _{n}\right) 
\notag \\
&&\ =\left( q,q^{1/2}e^{2i\theta },\ q^{1/2}e^{-2i\theta };\ q\right)
_{\infty },  \notag
\end{eqnarray}
for $\mathcal{E}_{q}\left( x;i\omega _{n}\right) $. If $\varphi =\pi /2$,
one gets 
\begin{eqnarray}
&&\sum_{n=-\infty }^{\infty }\frac{\left( -qr^{2}\omega
_{n}^{2};q^{2}\right) _{\infty }}{\left( -q\omega _{n}^{2};q^{2}\right)
_{\infty }}\ k^{-1}\left( \omega _{n}\right) \\
&&\ \qquad \times \ \mathcal{E}_{q}\left( \cos \theta ;i\omega _{n}\right) 
\notag \\
&&\ =\frac{\left( q,r^{2},q^{1/2}e^{2i\theta },\ q^{1/2}e^{-2i\theta };\
q\right) _{\infty }}{\pi \left( -r^{2}e^{2i\theta },-r^{2}e^{-2i\theta
};q^{2}\right) _{\infty }}.  \notag
\end{eqnarray}
A terminating case of this generating relation appears when $%
r^{2}=-1/q\omega _{m}^{2}$ for an integer $m\neq 0$, 
\begin{eqnarray}
&&\sum_{n=-\left| m\right| }^{\left| m\right| }\frac{\left( \omega
_{n}^{2}/\omega _{m}^{2};q^{2}\right) _{\infty }}{\left( -q\omega
_{n}^{2};q^{2}\right) _{\infty }}\ k^{-1}\left( \omega _{n}\right) \\
&&\ \qquad \times \ \mathcal{E}_{q}\left( \cos \theta ;i\omega _{n}\right) 
\notag \\
&&\ =\frac{\left( q,r^{2},q^{1/2}e^{2i\theta },\ q^{1/2}e^{-2i\theta };\
q\right) _{\infty }}{\pi \left( e^{2i\theta }/q\omega _{m}^{2},e^{-2i\theta
}/q\omega _{m}^{2};q^{2}\right) _{\infty }}.  \notag
\end{eqnarray}
Here $m=\pm 1,\pm 2,\pm 3,...$.

\section{Method of Summation of Basic Fourier Series}

According to Theorem 9.5, for a continuous function $f(x)$ the basic Fourier
series $\mathbf{S}\left[ f\right] $ converges to $f(x)$ if it converges
uniformly. In this section we shall discuss another method of summation of
basic Fourier series.

Let $f(x)$ be a bounded function that is continuous on $\left( -1,1\right) $
and let $\mathbf{S}\left[ f\right] $ be its $q$-Fourier series defined by
the right side of (8.5). Replace this series by 
\begin{equation}
\mathbf{S}_{r}\left[ f\right] =\sum_{n=-\infty }^{\infty }\ c_{n}(r)\ 
\mathcal{E}_{q}\left( \cos \theta ;i\omega _{n}\right) ,
\end{equation}
where 
\begin{eqnarray}
c_{n}(r) &=&\frac{\left( -qr^{2}\omega _{n}^{2};q^{2}\right) _{\infty }}{%
\left( -q\omega _{n}^{2};q^{2}\right) _{\infty }}\  \\
&&\times \ \frac{1}{2k\left( \omega _{n}\right) }\ \int_{0}^{\pi }f\left(
\cos \theta \right) \ \mathcal{E}_{q}\left( \cos \theta ;-ir\omega
_{n}\right)  \notag \\
&&\qquad \qquad \qquad \times \ \frac{\ \left( e^{2i\theta },\ e^{-2i\theta
};\ q\right) _{\infty }}{\left( q^{1/2}e^{2i\theta },\ q^{1/2}e^{-2i\theta
};\ q\right) _{\infty }}\ d\theta  \notag
\end{eqnarray}
provided that $0<r<1$. Comparing (11.2) and (8.6), 
\begin{equation}
\lim_{r\rightarrow 1^{-}}c_{n}(r)=c_{n},
\end{equation}
where $c_{n}$ are the regular $q$-Fourier coefficients of $f(x)$. Suppose
that the series $\mathbf{S}_{r}\left[ f\right] $ converges uniformly with
respect to the parameter $r$ when $0<r<1$. Then, 
\begin{equation}
\lim_{r\rightarrow 1^{-}}\mathbf{S}_{r}\left[ f\right] =\mathbf{S}\left[
f\right] .
\end{equation}

On the other hand, from (11.1)--(11.2) one gets 
\begin{eqnarray}
\mathbf{S}_{r}\left[ f\right] &=&\sum_{n=-\infty }^{\infty }\frac{\left(
-qr^{2}\omega _{n}^{2};q^{2}\right) _{\infty }}{\left( -q\omega
_{n}^{2};q^{2}\right) _{\infty }}\ \mathcal{E}_{q}\left( \cos \theta
;i\omega _{n}\right) \\
&&\times \frac{1}{2k\left( \omega _{n}\right) }\ \int_{0}^{\pi }f\left( \cos
\varphi \right) \ \mathcal{E}_{q}\left( \cos \varphi ;-ir\omega _{n}\right) 
\notag \\
&&\qquad \quad \quad \qquad \times \ \frac{\ \left( e^{2i\varphi },\
e^{-2i\varphi };\ q\right) _{\infty }}{\left( q^{1/2}e^{2i\varphi },\
q^{1/2}e^{-2i\varphi };\ q\right) _{\infty }}\ d\varphi .  \notag
\end{eqnarray}
Using the uniform convergence of the series in the bilinear generating
function (10.1), we finally obtain 
\begin{equation}
\mathbf{S}_{r}\left[ f\right] =\frac{1}{2\pi }\int_{0}^{\pi }\frac{f\left(
\cos \varphi \right) \ \left( q,r^{2},e^{2i\varphi },\ e^{-2i\varphi };\
q\right) _{\infty }}{\left( re^{i\theta +i\varphi },re^{i\theta -i\varphi
},re^{i\varphi -i\theta },re^{-i\theta -i\varphi };q\right) _{\infty }}\
d\varphi .
\end{equation}
It has been shown in \cite{As:Ra:Su} (see also \cite{Wi}) that 
\begin{equation}
\lim_{r\rightarrow 1^{-}}\frac{1}{2\pi }\int_{0}^{\pi }\frac{f\left( \cos
\varphi \right) \ \left( q,r^{2},e^{2i\varphi },\ e^{-2i\varphi };\ q\right)
_{\infty }}{\left( re^{i\theta +i\varphi },re^{i\theta -i\varphi
},re^{i\varphi -i\theta },re^{-i\theta -i\varphi };q\right) _{\infty }}\
d\varphi =f\left( \cos \theta \right)
\end{equation}
for every bounded function $f\left( \cos \theta \right) $ that is continuous
on $0<\theta <\pi $. As a result we have proved the following theorem.

\begin{theorem}
Let $f(x)$ be a bounded function that is continuous on $\left( -1,1\right) $
and let $\mathbf{S}_{r}\left[ f\right] $ be the series defined by
(11.1)--(11.2). If $\mathbf{S}_{r}\left[ f\right] $ converges uniformly with
respect to the parameter $r$ when $0<r<1$, then $\lim_{r\rightarrow 1^{-}}%
\mathbf{S}_{r}\left[ f\right] =\mathbf{S}\left[ f\right] =f(x)$.
\end{theorem}

\section{Relation Between $q$-Trigonometric System and $q$-Legendre
Polynomials}

The trigonometric system $\left\{ e^{i\pi nx}\right\} _{n=-\infty }^{\infty
} $ and the system of the Legendre polynomials $\left\{ P_{m}\left( x\right)
\right\} _{m=0}^{\infty }$ are two complete systems in $L^{2}\left(
-1,1\right) $. The corresponding unitary transformation between these two
orthogonal basises and its inverse are 
\begin{equation}
e^{i\pi nx}=\left( \frac{2}{\pi n}\right) ^{1/2}\sum_{m=0}^{\infty
}i^{m}\left( m+1/2\right) \ J_{m+1/2}\left( \pi n\right) \ P_{m}\left(
x\right)
\end{equation}
and 
\begin{equation}
P_{m}\left( x\right) =\sum_{n=-\infty }^{\infty }\left( -i\right) ^{m}\left( 
\frac{1}{2\pi n}\right) ^{1/2}\ J_{m+1/2}\left( \pi n\right) \ e^{i\pi nx},
\end{equation}
respectively. Relation (12.1) is a special case of a more general expansion, 
\begin{equation}
e^{irx}=\left( \frac{2}{r}\right) ^{\nu }\Gamma \left( \nu \right)
\sum_{m=0}^{\infty }i^{m}\left( \nu +m\right) \ J_{\nu +m}\left( r\right) \
C_{m}^{\nu }\left( x\right) ,
\end{equation}
where $C_{m}^{\nu }\left( x\right) $ are ultraspherical polynomials and $%
J_{\nu +m}\left( r\right) $ are Bessel functions \cite{Wa}. Expansion (12.2)
is the Fourier series of the Legendre polynomials on $\left( -1,1\right) $.
Orthogonality properties of the trigonometric system and Legendre
polynomials lead to the orthogonality relations, 
\begin{equation}
\sum_{m=0}^{\infty }\frac{m+1/2}{\pi n}\ J_{m+1/2}\left( \pi n\right) \
J_{m+1/2}\left( \pi l\right) =\delta _{nl},
\end{equation}
\begin{equation}
\sum_{n=-\infty }^{\infty }\frac{m+1/2}{\pi n}\ J_{m+1/2}\left( \pi n\right)
\ J_{p+1/2}\left( \pi n\right) =\delta _{mp},
\end{equation}
for the corresponding Bessel functions.

The basic trigonometric system $\left\{ \mathcal{E}_{q}\left( x;i\omega
_{n}\right) \right\} _{n=-\infty }^{\infty }$ and the system of the
continuous $q$-ultraspherical polynomials $\left\{ C_{m}\left( \left.
x;\beta \right| q\right) \right\} _{m=0}^{\infty }$ with $\beta =q^{1/2}$,
which are the basic analogs of the Legendre polynomials, are two complete
orthogonal systems in $L_{\rho }^{2}\left( -1,1\right) $, where $\rho $ is
the weight function in the orthogonality relation (7.5). Therefore, there
exists a $q$-version of the unitary transformation (12.1)--(12.2).

Ismail and Zhang \cite{Is:Zh} have found the following $q$-analog of (12.3), 
\begin{eqnarray}
\mathcal{E}_{q}\left( x;i\omega \right) &=&\frac{\left( q;q\right) _{\infty
}\omega ^{-\nu }}{\left( -q\omega ^{2};q^{2}\right) _{\infty }\left( q^{\nu
};q\right) _{\infty }}   \notag\\
&&\times \sum_{m=0}^{\infty }i^{m}\left( 1-q^{\nu +m}\right) \ q^{m^{2}/4}\
J_{\nu +m}^{\left( 2\right) }\left( 2\omega ;q\right) \ C_{m}\left( \left.
x;q^{\nu }\right| q\right) ,
\end{eqnarray}
where $J_{\nu +m}^{\left( 2\right) }\left( 2\omega ;q\right) $ is Jackson's $%
q$-Bessel function (see, for example, \cite{Ga:Ra}). Special case $\nu =1/2$
gives the basic analog of the expansion (12.1), 
\begin{eqnarray}
\mathcal{E}_{q}\left( x;i\omega _{n}\right) &=&\frac{\left( q;q\right)
_{\infty }\omega _{n}^{-1/2}}{\left( -q\omega _{n}^{2};q^{2}\right) _{\infty
}\left( q^{1/2};q\right) _{\infty }}  \notag \\
&&\times \sum_{m=0}^{\infty }i^{m}\left( 1-q^{m+1/2}\right) \ q^{m^{2}/4}\
J_{m+1/2}^{\left( 2\right) }\left( 2\omega _{n};q\right) \ C_{m}\left(
\left. x;q^{1/2}\right| q\right) ,
\end{eqnarray}
where $\omega _{-n}=-\omega _{n}$ and $\omega _{0}=0,\omega _{1},\omega
_{2},\omega _{3},...$, are nonnegative zeros of the basic sine function $%
S\left( \eta ;\omega \right) $ arranged in ascending order of magnitude.

On the other hand, the continuous $q$-ultraspherical polynomials $%
C_{m}\left( \left. x;q^{1/2}\right| q\right) $ can be expanded in the $q$%
-Fourier series as 
\begin{eqnarray}
C_{m}\left( \left. x;q^{1/2}\right| q\right) &=&\pi \frac{\left(
q^{1/2};q\right) _{\infty }}{\left( q;q\right) _{\infty }} \notag \\
&&\times \sum_{n=-\infty }^{\infty }\left( -i\right) ^{m}q^{m^{2}/4}\ \frac{%
\omega _{n}^{-1/2}}{k\left( \omega _{n}\right) \left( -q\omega
_{n}^{2};q^{2}\right) _{\infty }}\ \ J_{m+1/2}^{\left( 2\right) }\left(
2\omega _{n};q\right)  \\
&&\qquad \qquad \qquad \qquad \times \ \mathcal{E}_{q}\left( x;i\omega
_{n}\right) .  \notag
\end{eqnarray}
Indeed, by (8.5)--(8.6), 
\begin{equation}
C_{m}\left( \left. x;q^{1/2}\right| q\right) =\sum_{n=-\infty }^{\infty
}c_{n}\ \mathcal{E}_{q}\left( x;i\omega _{n}\right) ,
\end{equation}
where 
\begin{eqnarray}
c_{n} &=&\frac{1}{2k\left( \omega _{n}\right) }\ \int_{0}^{\pi }C_{m}\left(
\left. \cos \theta ;q^{1/2}\right| q\right) \ \mathcal{E}_{q}\left( \cos
\theta ;-i\omega _{n}\right) \\
&&\qquad \quad \qquad \qquad \times \ \frac{\ \left( e^{2i\theta },\
e^{-2i\theta };\ q\right) _{\infty }}{\left( q^{1/2}e^{2i\theta },\
q^{1/2}e^{-2i\theta };\ q\right) _{\infty }}\ d\theta .  \notag
\end{eqnarray}
Using (12.7), where the series on the right converge uniformly in $x$ for
any $\omega $, and the orthogonality relation 
\begin{eqnarray}
&&\int_{0}^{\pi }C_{m}\left( \left. \cos \theta ;q^{1/2}\right| q\right) \
C_{p}\left( \left. \cos \theta ;q^{1/2}\right| q\right) \\
&&\qquad \quad \times \frac{\ \left( e^{2i\theta },\ e^{-2i\theta };\
q\right) _{\infty }}{\left( q^{1/2}e^{2i\theta },\ q^{1/2}e^{-2i\theta };\
q\right) _{\infty }}\ d\theta  \notag \\
&&\qquad =2\pi \frac{\left( q^{1/2};q\right) _{\infty }^{2}}{\left(
q;q\right) _{\infty }^{2}}\ \left( 1-q^{m+1/2}\right) ^{-1}\delta _{mp} 
\notag
\end{eqnarray}
(see, for example, \cite{Ga:Ra}), one gets 
\begin{eqnarray}
c_{n} &=&\pi \frac{\left( q^{1/2};q\right) _{\infty }}{\left( q;q\right)
_{\infty }}\ \left( -i\right) ^{m}q^{m^{2}/4} \\
&&\times \ \frac{\omega _{n}^{-1/2}}{k(\omega _{n})\left( -q\omega
_{n}^{2};q^{2}\right) _{\infty }}\ J_{m+1/2}^{\left( 2\right) }\left(
2\omega _{n};q\right) ,  \notag
\end{eqnarray}
or, 
\begin{eqnarray}
c_{n} &=&\pi \frac{\left( q^{1/2};q\right) _{\infty }^{2}}{\left( q;q\right)
_{\infty }^{2}}\ \frac{\left( -i\right) ^{m}q^{m^{2}/4}}{\left(
q^{1/2};q\right) _{m+1}}\ \frac{\ \omega _{n}^{m}\left( -\omega
_{n}^{2};q^{2}\right) _{\infty }}{k\left( \omega _{n}\right) \left( -q\omega
_{n}^{2};q^{2}\right) _{\infty }}  \notag \\
&&\times \ _{2}\varphi _{1}\left( 
\begin{array}{c}
\begin{array}{ll}
-q^{m+3/2}, & -q^{m+5/2}
\end{array}
\\ 
q^{2m+3}
\end{array}
;\ q^{2},\ -\omega _{n}^{2}\right) \\
&=&\pi \frac{\left( q^{1/2};q\right) _{\infty }^{2}}{\left( q;q\right)
_{\infty }^{2}}\ \frac{\left( -i\right) ^{m}q^{m^{2}/4}\ \omega _{n}^{m}}{%
\left( q^{1/2};q\right) _{m+1}\ k\left( \omega _{n}\right) }  \notag \\
&&\times \Biggl[\frac{\left( -q^{m+3/2};q\right) _{\infty }}{\left(
q,q^{2m+3};q^{2}\right) _{\infty }}\ \frac{\left( \omega
_{n}^{2}q^{m+3/2},q^{1/2-m}/\omega _{n}^{2};q^{2}\right) _{\infty }}{\left(
-q\omega _{n}^{2},-q^{2}/\omega _{n}^{2};q^{2}\right) _{\infty }}  \notag \\
&&\qquad \quad \times \ _{2}\varphi _{1}\left( 
\begin{array}{c}
\begin{array}{ll}
-q^{m+3/2}, & -q^{1/2-m}
\end{array}
\\ 
q
\end{array}
;\ q^{2},\ -\frac{q}{\omega _{n}^{2}}\right)  \notag \\
&&\qquad +\frac{\left( -q^{m+1/2};q\right) _{\infty }}{\left(
q^{-1},q^{2m+3};q^{2}\right) _{\infty }}\ \frac{\left( \omega
_{n}^{2}q^{m+5/2},q^{-1/2-m}/\omega _{n}^{2};q^{2}\right) _{\infty }}{\left(
-q\omega _{n}^{2},-q^{2}/\omega _{n}^{2};q^{2}\right) _{\infty }}  \notag \\
&&\qquad \quad \times \ _{2}\varphi _{1}\left( 
\begin{array}{c}
\begin{array}{ll}
-q^{m+5/2}, & -q^{3/2-m}
\end{array}
\\ 
q^{3}
\end{array}
;\ q^{2},\ -\frac{q}{\omega _{n}^{2}}\right) \Biggr],
\end{eqnarray}
by (5.31) and (III.32) of \cite{Ga:Ra}, respectively. The last equation
gives the large $\omega $-asymptotic of the basic Fourier coefficients. With
the aid of (5.11), (6.9), and (I.9) of \cite{Ga:Ra}, we finally obtain 
\begin{equation}
\left| c_{n}\right| \thicksim D\ q^{n/2}\rightarrow 0
\end{equation}
as $n\rightarrow \infty $, where $D$ is some constant. Therefore, the series
on the right side of (12.8) converges uniformly and we have established the
expansion of the $q$-Legendre polynomials $C_{m}\left( \left.
x;q^{1/2}\right| q\right) $ in terms of the basic trigonometric functions $%
\mathcal{E}_{q}\left( x;i\omega _{n}\right) $ due to Theorem 9.5.

Relations (12.7)--(12.8) define the unitary operator acting in $L_{\rho
}^{2}\left( -1,1\right) $ \cite{Akh:Gl}. Orthogonality relations of the
matrix of this operator lead to the following orthogonality properties 
\begin{eqnarray}
&&\sum_{m=0}^{\infty }\frac{\pi \left( 1-q^{m+1/2}\right) }{\omega _{n}\
k(\omega _{n})\ \left( -q\omega _{n}^{2};q^{2}\right) _{\infty }^{2}}\
q^{m^{2}/2} \\
&&\quad \times \ J_{m+1/2}^{\left( 2\right) }\left( 2\omega _{n};q\right) \
J_{m+1/2}^{\left( 2\right) }\left( 2\omega _{l};q\right) =\delta _{nl} 
\notag
\end{eqnarray}
and 
\begin{eqnarray}
&&\sum_{n=-\infty }^{\infty }\frac{\pi \left( 1-q^{m+1/2}\right) }{\omega
_{n}\ k(\omega _{n})\ \left( -q\omega _{n}^{2};q^{2}\right) _{\infty }^{2}}\
q^{m^{2}/2} \\
&&\quad \times \ J_{m+1/2}^{\left( 2\right) }\left( 2\omega _{n};q\right) \
J_{p+1/2}^{\left( 2\right) }\left( 2\omega _{n};q\right) =\delta _{mp} 
\notag
\end{eqnarray}
for the corresponding Jackson's $q$-Bessel function. These relations are,
clearly, $q$-analogs of (12.4)--(12.5).

\section{Some Basic Trigonometric Identities}

One of the most important formulas for the trigonometric functions is the
main trigonometric identity, 
\begin{equation}
\cos ^{2}\omega x+\sin ^{2}\omega x=1.
\end{equation}
It follows from the Pythagorean Theorem or from the addition formulas for
the trigonometric functions, but one can also prove this identity on the
base of the differential equation. The functions $\cos \omega x$ and $\sin
\omega x$ are two solutions of (1.9) corresponding to the same eigenvalue $%
\omega $. Therefore, 
\begin{equation}
\frac{d}{dx}\ \left[ W\left( \cos \omega x,\ \sin \omega x\right) \right] =0,
\end{equation}
or 
\begin{equation}
\cos ^{2}\omega x+\sin ^{2}\omega x=\text{constant}.
\end{equation}
Substituting $x=0$, one verifies (13.1).

One can extend this consideration to the case of the basic trigonometric
functions. Consider equation (3.8) with $u(z)=C_{q}\left( x(z);\omega
\right) $, $v(z)=S_{q}\left( x(z);\omega \right) $, and $\rho \left(
z\right) =1$, 
\begin{equation}
\Delta \left[ W\left( u(z),\ v(z)\right) \right] =0,
\end{equation}
where 
\begin{eqnarray}
W\left( u,\ v\right) &=&W\left( C\left( x;\omega \right) ,\ S\left( x;\omega
\right) \right) \\
&=&\frac{2q^{1/4}\omega }{1-q}\ \left[ C\left( x(z);\omega \right) C\left(
x(z-1/2);\omega \right) \right.  \notag \\
&&\qquad \ \ \ +\left. S\left( x(z);\omega \right) S\left( x(z-1/2);\omega
\right) \right]  \notag
\end{eqnarray}
is the analog of the Wronskian (3.9) and we also used (2.12)--(2.13). One
can easily see that $W\left( u,\ v\right) $ here is a doubly periodic
function in $z$ without poles in the rectangle on the Figure. Therefore,
this function is just a constant by Liouville's theorem, 
\begin{equation*}
C\left( x\left( z\right) \right) C\left( x\left( z-1/2\right) \right)
+S\left( x\left( z\right) \right) S\left( x\left( z-1/2\right) \right) =C.
\end{equation*}
The value of this constant $C$ can be found by choosing $x=0$, which gives 
\begin{eqnarray*}
C &=&\frac{\left( -\omega ^{2};q^{2}\right) _{\infty }^{2}}{\left( -q\omega
^{2};q^{2}\right) _{\infty }^{2}} \\
&&\times \ _{2}\varphi _{1}\left( 
\begin{array}{c}
\begin{array}{ll}
q, & q
\end{array}
\\ 
q
\end{array}
;\ q^{2},\ -\omega ^{2}\right) \\
&&\times \ _{2}\varphi _{1}\left( 
\begin{array}{c}
\begin{array}{ll}
1, & q
\end{array}
\\ 
q
\end{array}
;\ q^{2},\ -\omega ^{2}\right) \\
&=&\frac{\left( -\omega ^{2};q^{2}\right) _{\infty }^{2}}{\left( -q\omega
^{2};q^{2}\right) _{\infty }^{2}} \\
&&\times \ _{1}\varphi _{0}\left( 
\begin{array}{c}
q \\ 
-
\end{array}
;\ q^{2},\ -\omega ^{2}\right) \\
&=&\frac{\left( -\omega ^{2};q^{2}\right) _{\infty }}{\left( -q\omega
^{2};q^{2}\right) _{\infty }}
\end{eqnarray*}
by the $q$-binomial theorem. As a result one gets 
\begin{eqnarray}
&&C_{q}\left( \cos \theta ;\omega \right) C_{q}\left( \cos \left( \theta
+i\log q/2\right) ;\omega \right) \\
&&+\ S_{q}\left( \cos \theta ;\omega \right) S_{q}\left( \cos \left( \theta
+i\log q/2\right) ;\omega \right) =\frac{\left( -\omega ^{2};q^{2}\right)
_{\infty }}{\left( -q\omega ^{2};q^{2}\right) _{\infty }}  \notag
\end{eqnarray}
as a $q$-extension of the main identity (13.1). The special case $z=1/4$,
when $\eta =x(1/4)$, of (13.6) has the simplest form 
\begin{equation}
\ C^{2}\left( \eta ;\omega \right) +S^{2}\left( \eta ;\omega \right) =\frac{%
\left( -\omega ^{2};q^{2}\right) _{\infty }}{\left( -q\omega
^{2};q^{2}\right) _{\infty }}.
\end{equation}
Our identity (13.6) can also be derived as a special case of the
``addition'' theorem for the basic trigonometric functions established in 
\cite{Su}.

In a similar fashion, we can find an analog of the identity 
\begin{equation}
\cos ^{2}\omega \left( x+y\right) +\sin ^{2}\omega \left( x+y\right) =1
\end{equation}
considering more general basic sine and cosine functions, $C\left(
x,y;\omega \right) $ and $S\left( x,y;\omega \right) $, as two solutions of
equation (2.10). The result is 
\begin{eqnarray}
&&C_{q}\left( \cos \theta ,\cos \varphi ;\omega \right) C_{q}\left( \cos
\left( \theta +i\log q/2\right) ,\cos \varphi ;\omega \right) \\
&&+\ S_{q}\left( \cos \theta ,\cos \varphi ;\omega \right) S_{q}\left( \cos
\left( \theta +i\log q/2\right) ,\cos \varphi ;\omega \right)  \notag \\
\ \ &&\ =\frac{\left( -\omega ^{2};q^{2}\right) _{\infty }}{\left( -q\omega
^{2};q^{2}\right) _{\infty }}\ \left[ C_{q}^{2}\left( \cos \varphi ;\omega
\right) +S_{q}^{2}\left( \cos \varphi ;\omega \right) \right]  \notag \\
\ &&\ =\frac{\left( -\omega ^{2};q^{2}\right) _{\infty }}{\left( -q\omega
^{2};q^{2}\right) _{\infty }}\ C_{q}\left( \cos \varphi ,-\cos \varphi
;\omega \right) .  \notag
\end{eqnarray}
We have used (6.2) here. This identity can also be verified with the aid of
the ``addition'' theorems for the basic trigonometric functions.

Identity (13.7) gives the values of the basic cosine function $C\left( \eta
;\omega \right) $ at the zeros of the basic sine function $S\left( \eta
;\omega \right) $, 
\begin{equation}
C\left( \eta ;\omega _{n}\right) =\left( -1\right) ^{n}\ \sqrt{\frac{\left(
-\omega _{n}^{2};q^{2}\right) _{\infty }}{\left( -q\omega
_{n}^{2};q^{2}\right) _{\infty }}},
\end{equation}
and vice versa, 
\begin{equation}
S\left( \eta ;\varpi _{n}\right) =\left( -1\right) ^{n}\ \sqrt{\frac{\left(
-\varpi _{n}^{2};q^{2}\right) _{\infty }}{\left( -q\varpi
_{n}^{2};q^{2}\right) _{\infty }}},
\end{equation}
with the aid of Theorem 5.4.

\section{Example}

Let us consider a periodic function $p_{1}(x)$ which is defined in the
interval $\left( -1,1\right) $ by $p_{1}(x)=x$. Its Fourier coefficients are 
\begin{eqnarray*}
c_{0} &=&0; \\
c_{n} &=&\frac{1}{2}\int_{-1}^{1}xe^{-i\pi nx}\ dx \\
&=&\frac{\left( -1\right) ^{n-1}}{i\pi n},\qquad n\neq 0.
\end{eqnarray*}
Therefore, 
\begin{eqnarray}
x &=&\left. \sum_{n=-\infty }^{\infty }\right. ^{\prime }\frac{\left(
-1\right) ^{n-1}}{2i\pi n}\ e^{i\pi nx} \\
&=&2\sum_{n=1}^{\infty }\left( -1\right) ^{n-1}\ \frac{\sin \pi nx}{\pi n}. 
\notag
\end{eqnarray}

The special case $m=1$ of (12.8), 
\begin{eqnarray}
C_{1}\left( \left. x;q^{1/2}\right| q\right) &=&-i\pi \frac{\left(
q^{1/2};q\right) _{\infty }}{\left( q;q\right) _{\infty }}\ q^{1/4} \\
&&\times \sum_{n=-\infty }^{\infty }\frac{\omega _{n}^{-1/2}}{k\left( \omega
_{n}\right) \left( -q\omega _{n}^{2};q^{2}\right) _{\infty }}\ \
J_{3/2}^{\left( 2\right) }\left( 2\omega _{n};q\right)  \notag \\
&&\qquad \qquad \qquad \qquad \times \ \mathcal{E}_{q}\left( x;i\omega
_{n}\right) ,  \notag
\end{eqnarray}
gives us a possibility to establish the $q$-analog of (14.1). Let us first
simplify the right side of (14.2). Using the three-term recurrence relation
for the $q$-Bessel functions (see Exercise 1.25 of \cite{Ga:Ra}) and
(5.32)--(5.33), one gets 
\begin{eqnarray}
J_{3/2}^{\left( 2\right) }\left( 2\omega _{n};q\right)
&=&-q^{-1/2}J_{-1/2}^{\left( 2\right) }\left( 2\omega _{n};q\right) \\
&=&-\frac{\left( q^{1/2};q\right) _{\infty }}{\left( q;q\right) _{\infty }}\ 
\frac{\left( -q\omega _{n}^{2};q^{2}\right) _{\infty }}{\left( q\omega
_{n}\right) ^{1/2}}\ C\left( \eta ;\omega _{n}\right) .  \notag
\end{eqnarray}
On the other hand, 
\begin{equation}
C_{1}\left( \left. x;q^{1/2}\right| q\right) =\frac{2}{1+q^{1/2}}\ x.
\end{equation}
Combining (14.2)--(14.4) and (13.10), we finally obtain, 
\begin{eqnarray}
x &=&\pi \frac{\left( q^{1/2};q\right) _{\infty }^{2}}{\left( q;q\right)
_{\infty }^{2}}\ \frac{1}{2}\left( q^{1/4}+q^{-1/4}\right) \\
&&\times \left. \sum_{n=-\infty }^{\infty }\right. ^{\prime }\frac{\left(
-1\right) ^{n-1}}{i\ k\left( \omega _{n}\right) \ \omega _{n}}\ \sqrt{\frac{%
\left( -\omega _{n}^{2};q^{2}\right) _{\infty }}{\left( -q\omega
_{n}^{2};q^{2}\right) _{\infty }}}\   \notag \\
&&\qquad \qquad \qquad \quad \times \ \mathcal{E}_{q}\left( x;i\omega
_{n}\right)  \notag \\
&=&\pi \frac{\left( q^{1/2};q\right) _{\infty }^{2}}{\left( q;q\right)
_{\infty }^{2}}\ \left( q^{1/4}+q^{-1/4}\right)  \notag \\
&&\times \sum_{n=1}^{\infty }\frac{\left( -1\right) ^{n-1}}{k\left( \omega
_{n}\right) \ \omega _{n}}\ \sqrt{\frac{\left( -\omega _{n}^{2};q^{2}\right)
_{\infty }}{\left( -q\omega _{n}^{2};q^{2}\right) _{\infty }}}  \notag \\
&&\qquad \qquad \qquad \quad \times \ S_{q}\left( x;\omega _{n}\right) . 
\notag
\end{eqnarray}
These equations are, clearly, $q$-analogs of (14.1).

\section{Miscellaneous Results}

Under certain restrictions a function $f(z)$ analytic in the entire complex
plane and having zeros at the points $a_{1},a_{2},a_{3},...$ (these are the
only zeros of $f(z)$), where $\lim_{n\rightarrow \infty }\left| a_{n}\right| 
$ is infinite, can be represented as an infinite product 
\begin{equation}
f(z)=f(0)\ e^{zf^{\prime }(0)/f(0)}\prod_{n=1}^{\infty }\left( \left( 1-%
\frac{z}{a_{n}}\right) e^{z/a_{n}}\right)
\end{equation}
see, for example, \cite{Wh:Wa}, \cite{Lev}. Consider the entire function 
\begin{equation}
f(\omega )=\left( -q\omega ^{2};q^{2}\right) _{\infty }\frac{S(\eta ;\omega )%
}{\omega }
\end{equation}
which has simple real zeros at $\omega =\pm \omega _{n}$ by Theorems
5.1--5.3. In this case 
\begin{eqnarray*}
f(0) &=&\frac{1}{1-q^{1/2}}, \\
f(\omega ) &=&f(0)+\frac{1}{2}f^{\prime \prime }(0)\ \omega ^{2}+..., \\
f^{\prime }(0) &=&0
\end{eqnarray*}
and 
\begin{equation*}
f(\omega )=\frac{1}{1-q^{1/2}}\prod_{n=1}^{\infty }\left( \left( 1-\frac{%
\omega }{\omega _{n}}\right) e^{\omega /\omega _{n}}\right) \left( \left( 1+%
\frac{\omega }{\omega _{n}}\right) e^{-\omega /\omega _{n}}\right) .
\end{equation*}
As a result we arrive at the infinite product representation for the basic
sine function, 
\begin{eqnarray}
S(\eta ;\omega ) &=&\frac{1}{1-q^{1/2}}\ \frac{\omega }{\left( -q\omega
^{2};q^{2}\right) _{\infty }} \\
&&\times \prod_{n=1}^{\infty }\left( \left( 1-\frac{\omega }{\omega _{n}}%
\right) e^{\omega /\omega _{n}}\right) \left( \left( 1+\frac{\omega }{\omega
_{n}}\right) e^{-\omega /\omega _{n}}\right)  \notag \\
&=&\frac{1}{1-q^{1/2}}\ \frac{\omega }{\left( -q\omega ^{2};q^{2}\right)
_{\infty }}\prod_{n=1}^{\infty }\left( 1-\frac{\omega ^{2}}{\omega _{n}^{2}}%
\right) .  \notag
\end{eqnarray}
In a similar manner, one can obtain an infinite product representation for
the basic cosine function, 
\begin{eqnarray}
C(\eta ;\omega ) &=&\frac{1}{\left( -q\omega ^{2};q^{2}\right) _{\infty }} \\
&&\times \prod_{n=1}^{\infty }\left( \left( 1-\frac{\omega }{\varpi _{n}}%
\right) e^{\omega /\varpi _{n}}\right) \left( \left( 1+\frac{\omega }{\varpi
_{n}}\right) e^{-\omega /\varpi _{n}}\right)  \notag \\
&=&\frac{1}{\left( -q\omega ^{2};q^{2}\right) _{\infty }}\
\prod_{n=1}^{\infty }\left( 1-\frac{\omega ^{2}}{\varpi _{n}^{2}}\right) . 
\notag
\end{eqnarray}

Equations (13.10)--(13.11) and (15.3)--(15.4) result in the following
relations, 
\begin{equation}
\left( -1\right) ^{n}\sqrt{\left( -\omega _{m}^{2};q\right) _{\infty }}%
=\prod_{n=1}^{\infty }\left( 1-\frac{\omega _{m}^{2}}{\varpi _{n}^{2}}\right)
\end{equation}
and 
\begin{equation}
\left( -1\right) ^{n}\sqrt{\left( -\varpi _{m}^{2};q\right) _{\infty }}%
=\prod_{n=1}^{\infty }\left( 1-\frac{\varpi _{m}^{2}}{\omega _{n}^{2}}%
\right) ,
\end{equation}
between the zeros of the basic sine $S(\eta ;\omega )$ and basic cosine $%
C(\eta ;\omega )$ functions.

\section{Appendix: Estimate of Number of Zeros of $S(\eta ;\omega )$}

In this section we give an estimate for number of zeros of the basic sine
function $S(\eta ;\omega )$ on the basis of Jensen's theorem (see, for
example, \cite{Boa} and \cite{Lev}). We shall apply the method proposed by
Mourad Ismail at the level of the third Jackson $q$-Bessel functions \cite
{Is}(see also \cite{Is:Ma:Su2} for an extension of his idea to $q$-Bessel
functions on a $q$-quadratic grid).

Let us consider the entire function $f(\omega )$ defined in (15.2) again and
let $n_{f}\left( r\right) $ be the number of of zeros of $f(\omega )$ in the
circle $\left| \omega \right| <r$. Consider also circles of radius $%
R=R_{n}=\varkappa q^{-n}$, $q^{1/4}\leq \varkappa <q^{-3/4}$ with $%
n=1,2,3,...$ in the complex $\omega $-plane.

Since $n_{f}\left( r\right) $ is nondecreasing with $r$ one can write 
\begin{equation}
n_{f}\left( R_{n}\right) \leq n_{f}\left( r\right) \leq n_{f}\left(
R_{n+1}\right)
\end{equation}
if $R_{n}\leq r\leq R_{n+1}$, and, therefore, 
\begin{equation}
n_{f}\left( R_{n}\right) \ \int_{R_{n}}^{R_{n+1}}\frac{dr}{r}\leq
\int_{R_{n}}^{R_{n+1}}\frac{n_{f}\left( r\right) }{r}\ dr\leq n_{f}\left(
R_{n+1}\right) \ \int_{R_{n}}^{R_{n+1}}\frac{dr}{r}.
\end{equation}
But 
\begin{equation*}
\left. \int_{R_{n}}^{R_{n+1}}\frac{dr}{r}=\log r\right|
_{R_{n}}^{R_{n+1}}=\log q^{-1}
\end{equation*}
and, finally, one gets 
\begin{equation}
\log q^{-1}\ n_{f}\left( R_{n}\right) \leq \int_{R_{n}}^{R_{n+1}}\frac{%
n_{f}\left( r\right) }{r}\ dr\leq \log q^{-1}\ n_{f}\left( R_{n+1}\right) .
\end{equation}
In the proof of Theorem 5.1 we have established the fact that for
sufficiently large $n$ there are at least two roots of $f(\omega )$ between
the circles $\left| \omega \right| =R_{n}$ and $\left| \omega \right|
=R_{n+1}$. Thus, for sufficiently large $n$ the inequality (16.3) should
really have one of the following forms 
\begin{equation}
\log q^{-1}\ n_{f}\left( R_{n}\right) \leq \int_{R_{n}}^{R_{n+1}}\frac{%
n_{f}\left( r\right) }{r}\ dr<\log q^{-1}\ n_{f}\left( R_{n+1}\right) ,
\end{equation}
or 
\begin{equation}
\log q^{-1}\ n_{f}\left( R_{n}\right) <\int_{R_{n}}^{R_{n+1}}\frac{%
n_{f}\left( r\right) }{r}\ dr\leq \log q^{-1}\ n_{f}\left( R_{n+1}\right) .
\end{equation}

Our next step is to estimate the integral in (16.4)--(16.5). By Jensen's
theorem \cite{Boa}, \cite{Lev} 
\begin{eqnarray}
\int_{R_{n}}^{R_{n+1}}\frac{n_{f}\left( r\right) }{r}\ dr
&=&\int_{0}^{R_{n+1}}\frac{n_{f}\left( r\right) }{r}\ dr-\int_{0}^{R_{n}}%
\frac{n_{f}\left( r\right) }{r}\ dr \\
&=&\frac{1}{2\pi }\int_{0}^{2\pi }\log \left| \frac{f\left( \varkappa
q^{-n-1}e^{i\vartheta }\right) }{f\left( \varkappa q^{-n}e^{i\vartheta
}\right) }\right| \ d\vartheta .  \notag
\end{eqnarray}
For large values of $n$ in view of (5.5), 
\begin{eqnarray*}
\frac{f\left( \varkappa q^{-n-1}e^{i\vartheta }\right) }{f\left( \varkappa
q^{-n}e^{i\vartheta }\right) } &\sim &\frac{\left( q^{3/2}\varkappa
^{2}q^{-2n-2}e^{2i\vartheta };q^{2}\right) _{\infty }}{\left(
q^{3/2}\varkappa ^{2}q^{-2n}e^{2i\vartheta };q^{2}\right) _{\infty }} \\
&=&1-q^{3/2}\varkappa ^{2}q^{-2n-2}e^{2i\vartheta },
\end{eqnarray*}
and 
\begin{equation*}
\log \left| \frac{f\left( \varkappa q^{-n-1}e^{i\vartheta }\right) }{f\left(
\varkappa q^{-n}e^{i\vartheta }\right) }\right| \sim 2n\ \log q^{-1}+\log
\alpha ,
\end{equation*}
where $\alpha =\varkappa ^{2}q^{-1/2}$. Therefore, 
\begin{equation}
\int_{R_{n}}^{R_{n+1}}\frac{n_{f}\left( r\right) }{r}\ dr=2n\ \log
q^{-1}+\log \alpha +\text{o}\left( 1\right)
\end{equation}
as $n\rightarrow \infty $.

From (16.3) and (16.7), 
\begin{equation*}
1+\frac{\log \alpha /\log q^{-1}}{2n}-\frac{1}{n}\leq \frac{n_{f}\left(
R_{n}\right) }{2n}\leq 1+\frac{\log \alpha /\log q^{-1}}{2n}
\end{equation*}
and, therefore, 
\begin{equation}
\lim_{n\rightarrow \infty }\frac{n_{f}\left( R_{n}\right) }{2n}=1.
\end{equation}
On the other hand, from (16.4)--(16.5), 
\begin{equation}
n_{f}\left( R_{n}\right) \leq 2n+\log \alpha /\log q^{-1}<n_{f}\left(
R_{n+1}\right)
\end{equation}
or 
\begin{equation}
n_{f}\left( R_{n}\right) <2n+\log \alpha /\log q^{-1}\leq n_{f}\left(
R_{n+1}\right) ,
\end{equation}
which gives 
\begin{eqnarray*}
&&n_{f}\left( R_{n+1}\right) -n_{f}\left( R_{n}\right) \\
&&\ <\left( 2n+2+\log \alpha /\log q^{-1}\right) -\left( 2n-2+\log \alpha
/\log q^{-1}\right) =4
\end{eqnarray*}
Thus, we have established that 
\begin{equation}
n_{f}\left( R_{n+1}\right) -n_{f}\left( R_{n}\right) <4
\end{equation}
as $n\rightarrow \infty $. Due to the symmetry $f\left( \omega \right) =$ $%
f\left( -\omega \right) $ the last inequality implies that there is only one
positive root of $S(\eta ;\omega )$ between the test points $\omega =\gamma
_{n\text{ }}$defined by (5.6) for large values of $n$.

\section{Acknowledgments}

We wish to thank Dick Askey, Mourad Ismail, John McDonald, and Mizan Rahman
for valuable discussions and comments. One of us (S. S.) gratefully
acknowledge the hospitality of the Department of Mathematics at Arizona
State University were this work was done.

\end{document}